\documentclass[11pt,a4paper]{article}
\usepackage{a4}
\usepackage{amssymb,latexsym}
\def\eqnreset{\setcounter{equation}{0}}
\def\eqsection#1{\section{#1}\eqnreset}

\newtheorem{Thm}{Theorem}[section]
\newtheorem{Defi}[Thm]{Definition}
\newtheorem{Cor}[Thm]{Corollary}
\newtheorem{Lemma}[Thm]{Lemma}
\newtheorem{Prop}[Thm]{Proposition}
\newtheorem{Rem}[Thm]{Remark}
\newtheorem{Conj}[Thm]{Conjecture}
\newtheorem{Prelim}[Thm]{Preliminary}

\newenvironment{thm}[0]{\begin{Thm}\noindent}%
{\end{Thm}}
\newenvironment{defi}[0]{\begin{Defi}\noindent\rm}%
{\end{Defi}}
\newenvironment{cor}[0]{\begin{Cor}\noindent}%
{\end{Cor}}
\newenvironment{lemma}[0]{\begin{Lemma}\noindent}%
{\end{Lemma}}
\newenvironment{prop}[0]{\begin{Prop}\noindent}%
{\end{Prop}}
\newenvironment{rem}[0]{\begin{Rem}\noindent\rm}%
{\end{Rem}}
{\end{Conj}}
{\end{Prelim}}

\def\proof{\par\noindent{\it Proof.}{\ }{\ }}
\def\qed{~\hfill$\square$\medbreak}

\def\naam#1{\label{#1}}
\def\refer#1{\ref{#1}}

\def\bib#1{\cite{#1}}

\def\text#1{\;\;\;\;{\rm \hbox{#1}}\;\;\;\;}
\def\qquad{\quad\quad}

\def\itema{\vspace{-1mm}\item[{\rm (a)}]}
\def\itemb{\item[{\rm (b)}]}

\def\msy#1{{\mathbb #1}}
\def\C{{\msy C}}
\def\N{{\msy N}}
\def\Z{{\msy Z}}
\def\R{{\msy R}}

\def\ga{\alpha}

\def\gf{\varphi}

\def\gl{\lambda}

\def\gs{\sigma}

\def\gS{\Sigma}

\def\got#1{\mathfrak #1}
\def\fa{{\got a}}

\def\fg{{\got g}}
\def\fh{{\got h}}

\def\fk{{\got k}}

\def\fm{{\got m}}
\def\fn{{\got n}}

\def\fp{{\got p}}
\def\fq{{\got q}}

\def\to{\rightarrow}
\def\Re{{\rm Re}\,}

\def\inp#1#2{\langle#1\,,\,#2\rangle}

\def\End{{\rm End}}
\def\Hom{{\rm Hom}}

\def\ad{{\rm ad}}
\def\after{\,{\scriptstyle\circ}\,}

\def\pr{{\rm pr}}

\def\pijl#1{{\buildrel #1 \over \longrightarrow}}
\def\tr{{\rm tr}\,}

\def\iq{{\rm q}}
\def\iC{{\scriptscriptstyle \C}}

\def\cA{{\mathcal A}}

\def\cC{{\mathcal C}}

\def\cF{{\mathcal F}}
\def\cH{{\mathcal H}}

\def\cL{{\mathcal L}}
\def\cM{{\mathcal M}}
\def\cN{{\mathcal N}}
\def\cO{{\mathcal O}}
\def\cP{{\mathcal P}}

\def\cT{{\mathcal T}}
\def\cU{{\mathcal U}}

\def\cW{{\mathcal W}}

\mathcode`:="603A
\def\col{\,:\,}
\def\too{\longrightarrow}
\def\Ci{C^\infty}
\def\Cartan{\theta}
\def\Aut{{\rm Aut}}
\def\fap{\fa_0}
\def\Ap{A_0}
\def\minparabs{\cP_0}
\def\fapd{\fa_0^*}
\def\fapdc{\fa_{0\iC}^*}

\def\stM{{M_0}}
\def\cAtwo{\cA_2}
\def\Cci{C_c^\infty}
\def\udE{{}^u\! E^*}
\def\Vtau{V_\tau}
\def\ev{{\rm ev}}
\def\uFou{{}^u \cF}
\def\tayl{{\rm tayl}}
\def\supp{{\rm supp}\,}
\def\uAC{{}^u\!{\rm AC}}

\def\dotvar{\,\cdot\,}
\def\EPW{{\mathrm H}}
\def\EPR{{\mathrm H}_R}
\def\uPW{{}^u {\rm PW}}

\def\embeds{\hookrightarrow}
\def\uPWd{{}^u{\rm PW}^*}
\def\PW{{\rm PW}}
\def\dE{E^*}
\def\Fou{\cF}

\def\multCP{\cU_P}
\def\multCPt{\cU_P^t}
\def\Mer{\cM}
\def\Hyp{\cH}
\def\laur{{\rm laur}}
\def\Lau{\cL}
\def\AC{{\rm AC}}
\def\Aut{{\rm Aut}}

\def\faQp{\fa_Q^+}
\def\AQp{A_Q^+}
\def\ooC{{}^\circ\hspace{-0.6pt} C}
\def\Ccminf{C^{- \infty}_c}
\def\hol{{\rm hol}}
\def\uAChol{{}^u{\rm AC}_\hol}
\def\void{\varnothing}

\def\mspace{\mbox{\hspace{-0.5pt}}}
\def\spaceind{{}_{\scriptscriptstyle *}}
\def\spgs{\spaceind\sigma}
\def\spG{\spaceind\mspace G}
\def\spH{\spaceind\mspace H}
\def\spK{\spaceind\mspace K}
\def\spfaq{{\spaceind \mspace\fa_\iq}}

\def\spgl{{\spaceind \gl}}
\def\spP{{\spaceind \mspace P}}
\def\sppsi{{\spaceind\mspace \psi}}
\def\nE{{E^\circ}}
\def\sptau{{\spaceind\mspace \tau}}
\def\spX{{\mathrm X}}
\def\spspX{{\spaceind \spX}}
\def\spf{{\spaceind f}}
\def\Cartan{\theta}
\def\spCartan{\spaceind\mspace\Cartan}
\def\minf{{-\infty}}
\def\spfp{\spaceind \fp}
\def\spfg{\spaceind \fg}
\def\spfq{\spaceind \fq}
\def\spfh{\spaceind \fh}

\def\barB{\bar B}
\def\spbarB{\spaceind \bar B}
\def\spfaqd{\spaceind \fa_\iq^*}
\def\spfaqdc{\spaceind \fa_{\iq\iC}^*}
\def\spga{\spaceind \ga}
\def\spgS{\spaceind \gS}
\def\spAq{\spaceind A_\iq}
\def\spPmin{\spaceind\cP_\gs^{\mathrm{min}}}
\def\sprho{\spaceind\rho}
\def\diag{\mathrm{diag}}
\def\spcAtwo{\spaceind\cA_2}
\def\spLau{\spaceind \cL}
\def\spW{\spaceind W}
\def\spcW{\spaceind \cW}
\def\oC{{}^\circ C}
\def\pspsi{\spaceind \psi}
\def\spM{\spaceind M}
\def\spbarP{\spaceind \bar P}
\def\tpsi{\tilde \psi}
\def\sptpsi{\spaceind \tilde \psi}
\def\HC{{\mbox{\tiny \rm HC}}}
\def\spHyp{\spaceind \cH}
\def\spd{\spaceind d}
\def\sppi{\spaceind \pi}
\def\fapdreg{\fa_0^{*{\rm reg}}}
\def\bfapd{\bar\fa_{0}^*}
\def\mult{\gamma}
\def\barB{\bar B}

\def\eFou{\cF_{\rm eucl}}
\begin{document}
\title{Paley--Wiener spaces for real reductive Lie groups}
\author{Erik P.\ van den Ban and Henrik Schlichtkrull}
\date{}
\maketitle
\noindent
{\small In honor of Gerrit van Dijk}
\tableofcontents
\eqsection{Introduction}
The Paley--Wiener theorem for $K$-finite
compactly supported smooth
functions
on a real reductive group Lie group $G$ of the Harish-Chandra class
is due to J. Arthur \bib{Arthur} in general, and to O.A.\ Campoli \bib{Camppw} for $G$ of split rank one.

In our paper \bib{BSpw} we established a similar Paley--Wiener theorem for
smooth functions on a reductive symmetric space. In this paper
we will show that Arthur's theorem is a consequence of our result if one
considers the group $G$ as a symmetric space for $G\times G$ with respect
to the left times right action. At the same time
we will formulate a Paley--Wiener theorem for $K$-finite
{\it generalized functions} (in the
sense of distribution theory) on $G$, and prove that it is a special case of the
Paley--Wiener theorem for symmetric spaces established  in our paper \bib{BSdpw}.

All mentioned Paley--Wiener theorems are formulated in the following spirit.
A Fourier transform is defined by means of Eisenstein integrals for the minimal
principal series for the group or space under consideration.
The Eisenstein integrals depend on a certain spectral parameter and
satisfy the so called Arthur--Campoli relations.
A Paley--Wiener space is defined as a certain space of meromorphic functions in the spectral parameter,
characterized by the Arthur--Campoli relations and by growth estimates.
The Paley--Wiener theorem
asserts that the Fourier transform is a topological linear isomorphism
from a space of $K$-finite compactly supported smooth or generalized functions onto
a particular Paley--Wiener space.

The Paley--Wiener space of Arthur's paper is defined in terms of Eisenstein
integrals as introduced by Harish-Chandra \bib{HC1}; we shall refer to these
integrals as being {\em unnormalized.}
Our  Paley--Wiener theorems
in \bib{BSpw} and \bib{BSdpw} are defined in terms of the so-called {\em normalized\ }
Eisenstein integrals defined in \bib{BSft}. For $G$ considered as a symmetric
space the normalized Eisenstein integral differs from the unnormalized one.
Consequently, the associated Fourier transforms and Paley--Wiener spaces
are different. The final objective of this paper
is to clarify
the relationships between the various  Paley--Wiener spaces.

Recently,
P.\ Delorme, \bib{Delpw}, has proved a different
Paley--Wiener theorem, involving the operator-valued
Fourier transforms associated  with  all  generalized
principal series representations. In his result the
Arthur--Campoli relations are replaced by intertwining relations.
Moreover, the result is valid without the restriction of $K$-finiteness.

We shall now give a brief outline of the present paper.
In Section \refer{s: basic notation} we introduce the basic concepts,
in particular the space $\Cci(G\col \tau)$
of $\tau$-spherical compactly supported smooth
functions,
for which Arthur's theorem is most
conveniently formulated.
We review the definition of Harish-Chandra's (unnormalized)
Eisenstein
integral $E(P\col \gl)$ for $P$ a minimal parabolic subgroup of $G$ and
with spectral parameter $\gl \in \fapdc;$ here $\fap$ is
the Lie algebra of the split component $A_0$ of $P.$ Finally, we give the definition
of the associated Fourier transform $\uFou_P.$

In Section \refer{s: Arthur's Paley--Wiener theorem}
we recall, in Theorem \refer{t: arthurs PW one}, the formulation
of Arthur's Paley--Wiener theorem, \bib{Arthur}.
This theorem deals with  the family of Fourier transforms $(\uFou_Q)_{Q \in \minparabs}$ with
$Q$ ranging over the finite set $\minparabs$ of parabolic subgroups with the same  split
component $\Ap.$ The theorem asserts that this  family of  transforms
establishes an isomorphism from $\Cci(G\col \tau)$
onto a Paley--Wiener space $\uPW(G,\tau, \minparabs).$
The Fourier transforms in the
family
are completely determined by
any single one among them. Accordingly, in Theorem \refer{t: arthurs PW two},
Arthur's Paley--Wiener theorem is
reformulated by asserting that a single Fourier transform $\uFou_P$ defines a topological
linear isomorphism from $\Cci(G\col \tau)$ onto a Paley--Wiener space
$\uPW(G,\tau, P).$

In Section \refer{s: dist PW space} we formulate, in Theorem \refer{t: distributional PW},
a distributional
Paley--Wiener theorem, which asserts that $\uFou_P$ defines a topological
linear isomorphism from the space $\Ccminf(G\col \tau)$ of $\tau$-spherical compactly supported
generalized functions on $G$ onto a Paley--Wiener space
$\uPW^*(G , \tau, P).$

At
this point of the paper, the main results have been stated. The rest
of the paper is devoted to proofs. In Section \refer{s: c functions} we
prepare for this by introducing the $C$-functions and listing
those of their properties  that are needed.

In Section \refer{s: Relations between the Fourier transforms}
we give the proof that  Theorem \refer{t: arthurs PW two} is indeed
a reformulation of Arthur's theorem.
We describe the relations between the different Fourier transforms
in terms of $C$-functions.
The equivalence of Theorems  \refer{t: arthurs PW one} and \refer{t: arthurs PW two}
is then captured by the assertion, in Proposition \refer{p: proj is topol iso},
that the natural map between the Paley--Wiener spaces $\uPW(G,\tau, \minparabs)$
and $\uPW(G, \tau, P)$ is a topological isomorphism.

In the next Section, \refer{s: the normalized Fourier transform},
a normalized Fourier transform $\Fou_P$ is defined in terms of the normalized Eisenstein integral
$$
\nE(P \col \gl \col \dotvar) := E(P \col \gl \col \dotvar) \after C_{P|P}(1 \col \gl)^{-1};
$$
This normalization is natural from the point of view of asymptotic expansions; it has the effect
that
the new, normalized $C$-functions become unitary for imaginary $\gl.$
It follows from the relation between the Eisenstein integrals that the normalized
Fourier transform is related to the above Fourier transform
by a relation of the form $\uFou_P = \multCP \after \Fou_P,$
where $\multCP$ denotes
multiplication by a  $C$-function
$\gl \mapsto C_{P|P}(1 \col - \bar \gl)^*.$
An associated Paley--Wiener space $\PW(G,\tau, P)$ is defined, as well as a distributional
Paley--Wiener space, indicated by superscript $*.$ The main result of the section, Theorem
\refer{t: multC is iso PW spaces}, asserts that the map $\multCP$ induces isomorphisms
$\PW(G,\tau, P) \to \uPW(G,\tau, P)$
and $\PW^*(G,\tau, P) \to \uPW^*(G,\tau, P).$
As a result, the associated Paley--Wiener theorems for the normalized Fourier transform are equivalent
to those for the unnormalized transform.

Let $\spspX$ denote $G$ viewed as a symmetric space for $\spG: = G \times G.$
In the final section the unnormalized Eisenstein integral $E(P\col \gl)$ for $G$
is related to the unnormalized Eisenstein integral for $\spspX$ as defined in \bib{BSft}.
It follows from this relation that the normalized Eisenstein integrals,
 for $G$ and $\spspX,$ coincide; therefore,
so do the normalized Fourier transforms. Likewise, it is shown that the Paley--Wiener spaces
for $G$ coincide with the similar spaces for $\spspX.$ This finally establishes the validity
of all mentioned Paley--Wiener theorems as special cases of the theorems in \bib{BSpw} and \bib{BSdpw}.

\eqsection{Basic concepts}
\naam{s: basic notation}
Let $G$ be a real reductive Lie group of the Harish-Chandra class
and let $K$ be a maximal compact subgroup. Let $\Vtau$ be a finite
dimensional Hilbert space, and let $\tau$ be a unitary
double representation
of $K$ in $\Vtau.$ By this we mean that $\tau = (\tau_1, \tau_2)$
with $\tau_1$ a left and $\tau_2$ a right unitary representation of $K$ in $\Vtau;$ moreover,
the representations $\tau_1$ and $\tau_2$ commute. We will often drop the subscripts
on $\tau_1$ and $\tau_2,$ writing
$$
\tau(k_1) v \tau(k_2) = \tau_1(k_1) v \tau_2(k_2)
$$
for all $v \in \Vtau$ and $k_1, k_2 \in K.$ A function $f: G \to \Vtau$ is called
$\tau$-spherical if it satisfies the rule
\begin{equation}
\naam{e: spherical rule}
f(k_1 g k_2) = \tau(k_1) f(g) \tau(k_2)
\end{equation}
for all $g \in G$ and $k_1, k_2 \in K.$ The space of smooth $\tau$-spherical functions
is denoted by $\Ci(G\col \tau)$ and equipped with the usual Fr\'echet topology.
The subspace $\Cci(G\col \tau)$ of compactly supported smooth $\tau$-spherical
functions is equipped with the usual complete locally convex (Hausdorff) topology.

We shall first briefly establish some notation for the group.
As usual, we denote Lie groups by Roman capitals, and their Lie algebras
by the Gothic lower case equivalents.

Let $\Cartan \in \Aut(G)$ be the Cartan involution associated with $K.$ The
associated infinitesimal involution of $\fg$ is denoted by the same symbol and
the associated eigenspaces with eigenvalues +1 and $-1$ by $\fk$ and $\fp,$ respectively.
Accordingly, we have the Cartan decomposition $\fg = \fk \oplus \fp.$

Let $\fap$ be a maximal abelian subspace of $\fp$ and let $\gS$
be the restricted root system of $\fap$ in $\fg.$
Let $\Ap:= \exp \fap$ be the associated
vectorial subgroup of $G$ and let $\minparabs =  \cP(\Ap)$ be the collection of parabolic subgroups
of $G$ with split component $\Ap.$ Each element $P \in \minparabs$
is minimal and has a Langlands decomposition of the form $P = \stM \Ap N_P,$
with $\stM$ equal to the centralizer of $\Ap$ in $K.$
Let $\gS(P)$ be the collection of $\fap$-roots in $\fn_P = {\rm Lie}(N_P).$
Then $P \mapsto \gS(P)$ defines a one-to-one correspondence between $\cP(A_0)$
and the collection of positive systems for $\gS.$

We equip $\fap$
with a $W$-invariant positive definite inner product $\inp{\dotvar}{\dotvar};$
the dual space $\fapd$ is equipped with the dual inner product. The latter inner product
is extended to a complex bilinear form  on $\fapdc.$ The norm associated with the inner product
on $\fapd$ is extended to a Hermitian
norm on $\fapdc,$ denoted by $|\dotvar|.$

We shall now review  the definition of the $\tau$-spherical Eisenstein
integral related to a given parabolic subgroup $P \in \cP(\Ap).$
Let $\tau_\stM$ denote the restriction of $\tau$ to $\stM.$ As $\stM$ is a subgroup of $K,$
the space $L^2(\stM \col \tau_\stM)$ of square integrable $\tau_\stM$-spherical
functions $\stM \to \Vtau $ is finite dimensional and equals the space of smooth $\tau_\stM$-spherical
functions. We equip $M_0$ with normalized Haar measure, and define the finite
dimensional Hilbert space
$\cAtwo = \cAtwo(\tau)$ by
\begin{equation}
\naam{e: defi cAtwo}
\cAtwo:= L^2(\stM\col \tau_\stM) =   \Ci(\stM\col \tau_\stM).
\end{equation}
Let $\psi \in \cAtwo.$  For $\gl \in \fapdc$
we define the function $\psi_\gl = \psi_{P, \gl}: G \to \Vtau$ by
$$
\psi_\gl(n a m k) = a^{\gl + \rho_P} \psi(m) \tau_2(k),
$$
for $k \in K,$ $m \in \stM,$ $a \in \Ap$ and $n \in N_P.$
Here $\rho_P \in \fapd$ is defined by
$$
\rho_P(H) = \frac12 \tr (\ad(H)|\fn_P).
$$
By the analytic nature of the Iwasawa decomposition $G = N_P \Ap K,$
the function $\psi_\gl$
is analytic. We define the Eisenstein integral
$
E(P\col \psi\col \gl) :\, G \to \Vtau
$
by
\begin{equation}
\naam{e: defi unnormalized E}
E(P\col \psi\col \gl \col x): = \int_K \tau_1(k)\psi_\gl(k^{-1} x) \; dk,
\end{equation}
for $x \in G.$ Then, clearly, $E(P\col \psi\col \gl)$ is a function in $\Ci(G\col \tau),$
depending linearly on $\psi$ and holomorphically on $\gl.$

\begin{rem}
\naam{r: gl parameter Eis}
Here we have adopted the same convention as J.\ Arthur \bib{Arthur}, \S 2, which differs from
Harish-Chandra's.
Let $E^\HC$ denote the Eisenstein integral as defined by Harish-Chandra \bib{HC1}, \S 9.
Then $E^\HC(P\col \psi \col \gl) = E(P \col \psi \col i \gl).$

For
 the reader's convenience, we note that Arthur  \bib{Arthur} uses the notation
$\cA_{\mathrm{cusp}}(M_0, \tau)$ or $\cA_0$ for the space (\refer{e: defi cAtwo}).
\end{rem}

In terms of the Eisenstein integral we define a Fourier transform
$\uFou_P$ from $\Cci(G\col \tau)$ to the space $\cO(\fapdc) \otimes \cAtwo$
of holomorphic $\cAtwo$-valued functions on $\fapdc.$ The superscript $u$ indicates that this Fourier
transform is defined in terms of the above unnormalized Eisenstein integral, in contrast with a normalized
Fourier transform $\Fou_P$ to be defined later.

Let $dx$ be a choice of Haar measure on
$G.$ We define the Fourier transform $\uFou_P f$ of  $f \in \Cci(G \col \tau)$ by the formula
\begin{equation}
\naam{e: defi unnormalized Fou}
\inp{\uFou_P f(\gl)}{\psi}_{\cAtwo} = \int_G \inp{f(x)}{E(P\col \psi\col - \bar\gl \col x)}_{\Vtau} \; dx,
\end{equation}
for $\psi \in \cAtwo$ and $\gl \in \fapdc.$ It follows from the Paley--Wiener theorem in \bib{Arthur} that
 $\uFou_P$ is injective on $\Cci(G\col \tau).$ This injectivity can also
be established by application
of the subrepresentation theorem, \bib{CM}, Thm.~8.21.

\eqsection{Arthur's Paley--Wiener theorem}
\naam{s: Arthur's Paley--Wiener theorem}
The image of $\Cci(G\col \tau)$ under Fourier transform
is described by the  Paley--Wiener theorem due to
J.\ Arthur \bib{Arthur}, which we shall now formulate.

It is convenient to rewrite the definition of the Fourier transform
$\uFou_P$ in terms of a (unnormalized) dual  Eisenstein integral.
Given $x \in G$ and $\gl \in \fapdc,$ we agree to define $E(P\col \gl\col x) \in \Hom(\cAtwo, \Vtau)$
by the formula
$$
E(P\col \gl \col x)\psi: = E(P \col \psi\col  \gl \col x)
$$
for $\psi \in \cAtwo.$ Moreover, we define a dual Eisenstein integral by
\begin{equation}
\naam{e: defi unnormalized dual Eis}
\udE(P \col \gl\col x): = E(P\col - \bar \gl \col x)^* \in \Hom(\Vtau, \cAtwo),
\end{equation}
where the superscript $*$ indicates that the Hilbert adjoint  has been taken.
The superscript $u$ serves to distinguish the present dual Eisenstein integral
from a normalized version that will be introduced at a later stage.

The dual Eisenstein integral $\udE$ may be viewed as a smooth
$\Hom(\Vtau, \cAtwo)$-valued
function on $\fapdc \times G$
which is holomorphic in the first variable.
In terms of this Eisenstein integral,
the Fourier transform
(\refer{e: defi unnormalized Fou})
may be expressed as an integral transform. Indeed,
it readily follows from the given definitions  that
\begin{equation}
\naam{e: defi uFou}
\uFou_P f (\gl) = \int_G \udE(P\col\gl \col x) f(x)\; dx,
\end{equation}
for $f \in \Cci(G\col \tau)$ and $\gl \in \fapdc.$
We now proceed to giving the definition of a suitable Paley--Wiener space.

Let $V$ be a finite dimensional real linear space.
We denote by $S(V)$ the symmetric algebra of $V_\iC.$ This algebra is identified
with the algebra of constant coefficient
holomorphic differential operators on $V_\iC$ in the usual way.

We  denote by $\cO(V_\iC)$
the space of holomorphic functions on $V_\iC$ and,
for $a \in V_\iC,$
by $\cO_a = \cO_a(V_\iC)$ the space of germs of holomorphic
functions at $a.$
Moreover, we denote by $\cO_a(V_\iC)^*_\tayl $ the space of
linear functionals
$\cO_a \to \C$ of the form
$$
f \mapsto \ev_a(uf) = u f(a),
$$
with $u \in S(V).$ The elements of $\cO_a(V_\iC)^*_\tayl$
will be called Taylor functionals
at $a,$ as they give  linear combinations
of coefficients of a Taylor series at $a.$
Clearly, the map $u \mapsto \ev_a\after u$ is a
linear isomorphism from $S(V)$ onto $\cO_a(V_\iC)^*_\tayl.$

We define the space of Taylor functionals on $V_\iC$
as the algebraic direct sum
$$
\cO(V_\iC)^*_\tayl \,:= \;\oplus_{a \in V_\iC}\;\; \cO_a(V_\iC)^*_\tayl.
$$
Given $U \in \cO(V_\iC)^*_\tayl,$ the finite set
of $a \in V_\iC$ with $U_a \neq 0$
is called the support of $U,$ notation $\supp U.$
Given $f \in \cO(V_\iC),$
we put
$$
U f := \sum_{a \in \supp U} U_a f_a.
$$
The map $U \otimes f \mapsto U f$ defines an embedding of
$\cO(V_\iC)^*_\tayl$
into the linear dual $\cO(V_\iC)^*;$
this justifies the notation.
Note that in fact the elements of $\cO(V_\iC)^*_\tayl$
are continuous with respect to
the usual Fr\'echet topology on $\cO(V_\iC).$

Finally, we note that a finitely supported function $U: V_\iC \to S(V)$ may be viewed
as a Taylor functional by the formula
$
U f := \sum_a \ev_a [ U(a) f ].
$
Accordingly, the space of Taylor functionals may be identified with
the space of finitely supported functions $V_\iC \to S(V).$

\begin{defi}
An (unnormalized, holomorphic)
Arthur--Campoli functional for $(G, \tau, \minparabs)$ is a family
$(\cL_P)_{P \in \minparabs} \subset \cO(\fapdc)^*_\tayl \otimes \cAtwo^*$ such that
$$
\sum_{P \in \minparabs} \cL_P\; [\udE (P\col \dotvar \col x)v_P] = 0
$$
for all $x \in G$ and all  $(v_P)_{P \in \minparabs} \subset \Vtau.$
The linear space of such families is denoted by
$\uAChol(G,\tau,\minparabs).$
\end{defi}

For $R > 0$ we define
$\EPR(\fapdc)$ to be the space of holomorphic functions
$\gf: \fapdc \to \C$ such that for every $n \in \N,$
$$
\nu_{R,n}(\gf) := \sup_{\gl \in \fapdc}\;\;
(1 + |\gl|)^{n} e^{- R |\Re \gl|}\, |\gf(\gl)| <  \infty.
$$
Equipped with the seminorms $\nu_{R, n},$ for $n \in \N,$
this space is a Fr\'echet space.

\begin{defi}
\naam{d: arthurs PW}
Let $R > 0.$ The Paley--Wiener space $\uPW_R(G, \tau, \minparabs)$ is defined to be the space of
families $(\gf_P)_{P \in \minparabs} \subset \EPR(\fapdc) \otimes \cAtwo$ such that
$$
\sum_{P \in \minparabs} \; \cL_P \gf_P = 0
$$
for all $(\cL_P)_{P \in \minparabs}  \in \uAChol(G,\tau, \minparabs).$
\end{defi}

By continuity of the Taylor functionals, the Paley--Wiener space is a closed subspace of the direct
sum of a finite number of copies of $ \EPR(\fapdc) \otimes \cAtwo,$ one for each $P \in \minparabs;$
it is therefore a Fr\'echet space of its own right.

For $R >0$ we put $G_R: = K \exp \barB_R  K,$ where $\barB_R$
 denotes the closed ball of center
$0$ and radius $R$ in $\fap.$ Moreover, we define the following closed subspace
of $\Cci(G\col \tau),$
and equip it with the relative topology:
$$
C_R^\infty(G \col \tau) = \{ f \in \Ci(G\col \tau) \mid \supp f \subset G_R \}.
$$
We have now gathered the concepts and notation needed to formulate
the Paley--Wiener theorem due to J.\ Arthur, \bib{Arthur}, p.~83, Thm.~3.3.l.

\begin{thm}{\rm (Arthur \bib{Arthur})\ }
\naam{t: arthurs PW one}
The map $f \mapsto (\uFou_P f)_{P \in \minparabs}$ is a topological linear isomorphism from
$\Ci_R(G\col \tau)$ onto $\uPW_R(G,\tau,\minparabs).$
\end{thm}

Each of the individual Fourier transforms $\uFou_P,$ for $P \in \minparabs,$
is already injective on $\Cci(G\col \tau).$ It is therefore natural
to reformulate Arthur's theorem in terms of a single Fourier transform.

\begin{defi}
Let $P \in \minparabs.$
An (unnormalized, holomorphic) Arthur--Campoli functional for the triple $(G,\tau, P)$
is a functional $\Lau \in \cO(\fapdc)^*_\tayl \otimes \cAtwo^*$ such that
\begin{equation}
\naam{e: Lau on udE}
\Lau [\udE(P \col \dotvar \col x)v] = 0,
\end{equation}
for all $x \in G$ and all $v \in \Vtau.$ The space of such functionals is denoted  by
$\uAChol(G,\tau, P).$
\end{defi}

\begin{defi}
\naam{d: unnormalized PW}
Let $P \in \minparabs$ and $R > 0.$ We define the
Paley--Wiener space $\uPW_R(G, \tau, P)$ to be the space
of functions $\gf \in \EPR(\fapdc) \otimes \cAtwo$ such that $\Lau \gf = 0$ for all
$\Lau \in \uAChol(G,\tau, P).$
\end{defi}

Arthur's Paley--Wiener theorem may now be reformulated as follows.

\begin{thm}
\naam{t: arthurs PW two}
Let $P \in \minparabs$ and $R > 0.$
The map $\uFou_P$ is a topological linear isomorphism
from $\Ci_R(G\col \tau)$ onto $\uPW_R(G,\tau,P).$
\end{thm}
The equivalence of Theorems
\refer{t: arthurs PW one} and \refer{t: arthurs PW two} will be established in
Section \refer{s: Relations between the Fourier transforms}.

\eqsection{A distributional Paley--Wiener space}
\naam{s: dist PW space}
In this section we will formulate a Paley--Wiener theorem
characterizing the image under Fourier transform
of the space $\Ccminf(G\col\tau)$ of compactly supported $\tau$-spherical
generalized  functions.

We shall first define the mentioned space.
The space $\Ccminf(G)$ of compactly supported
generalized functions on $G$ is defined as the topological linear dual
of the Fr\'echet space of smooth densities on $G.$
It is equipped with the strong dual topology.

Via the map $f \mapsto f \; dx,$ the space $\Ccminf(G)$ is isomorphic
with the space of compactly supported generalized densities on $G.$
Via integration the latter space may in turn be
identified with the continuous linear dual of $\Ci(G),$ i.e.,
with the space of compactly supported distributions on $G.$

The pairing with smooth densities induces
a natural embedding $\Cci(G) \embeds \Ccminf(G);$ accordingly,
the left and right regular representations of $G$ in $\Cci(G)$ extend
to  continuous representations of $G$ in $\Ccminf(G).$
We now
define $\Ccminf(G\col \Vtau)$
as the space of $K \times K$-invariants in $\Ccminf(G) \otimes \Vtau.$
Again, there is a natural embedding
$
\Cci(G\col \tau) \embeds \Ccminf(G \col \tau).
$

The definition
of $\uFou_P$ by (\refer{e: defi uFou}), for a given $P \in \minparabs,$
has a natural interpretation for
compactly supported $\tau$-spherical  generalized
 functions. Accordingly, $\uFou_P$ extends to a continuous
linear map $\Ccminf(G\col \tau) \to \cO(\fapdc) \otimes \cAtwo.$

Let $R > 0$ and $n \in \N.$ We define $\EPW^*_{R, n}(\fapdc)$ to
be the space of entire holomorphic functions $\gf: \fapdc \to \C$
with
$$
\nu_{R,- n}(\gf)= \sup_{\gl \in \fapdc}\; (1 + \|\gl\|)^{-n} e^{-R
| \Re \gl|} |\gf(\gl) |
 < \infty.
$$
Equipped with the given norm,  this space
is a Banach space. If $m < n,$ then
$$
\EPW^*_{R,m }(\fapdc) \subset  \EPW^*_{R, n}(\fapdc),
$$
with continuous inclusion map.
The union $\EPW_R^*(\fapdc)$ of these spaces,
for $n \in \N,$
is equipped with the inductive limit
locally convex topology.

\begin{defi}
\naam{d: unnormalized distrib PW}
Let $P \in \minparabs$ and $R > 0.$
The distributional Paley--Wiener space $\uPWd_R(G,\tau, P)$ is defined as the space
of functions $\gf \in \EPW_R^*(\fapdc) \otimes \cAtwo$ satisfying
the Arthur--Campoli relations $\Lau \gf = 0$ for all $\cL \in \uAChol(G,\tau, P).$
\end{defi}

By continuity of Taylor functionals, the space
$\uPWd_R(G, \tau, P)$ is a closed subspace
of $\EPW_R^*(\fapdc) \otimes \cAtwo;$ we endow it with the relative topology.

The following result is the analogue of Arthur's Paley--Wiener theorem
for $K$-finite compactly supported generalized functions on $G.$
Given $R > 0$ we denote by $C^{-\infty}_R(G\col \tau)$ the space
of $\tau$-spherical generalized functions on $G$ with support contained in $G_R = K\exp \barB_R K.$
It is equipped with the relative topology.

\begin{thm}
{\rm (The Paley--Wiener theorem for generalized functions)\ \ }
\naam{t: distributional PW}
Let $P \in \minparabs$ and $R > 0.$
The Fourier transform $\uFou_P$ extends to a topological
linear isomorphism from $C^{-\infty}_R(G \col \tau)$ onto $\uPW^*_R(G,\tau,P).$
\end{thm}

As mentioned in the introduction, this theorem will follow from the results of the present paper
combined with the distributional
Paley--Wiener theorem for  reductive symmetric
spaces proved in \bib{BSdpw}.

\begin{rem}
With the same arguments that will lead to the equivalence of Theorems \refer{t: arthurs PW two} and
\refer{t: arthurs PW one}, it can be shown that Theorem
\refer{t: distributional PW} is equivalent to a Paley--Wiener theorem for generalized functions
involving the family $(\uFou_P)_{P \in \minparabs}$ in the same spirit as Theorem \refer{t: arthurs PW one}.
\end{rem}

\eqsection{C-functions, singular loci and estimates}
\naam{s: c functions}
To establish the equivalence of Theorems \refer{t: arthurs PW one} and \refer{t: arthurs PW two}
we need relations between the Fourier transforms, which can be given in terms of the so-called $C$-functions.
The latter arise as coefficients in asymptotic expansions of Eisenstein integrals.

Let $Q \in \minparabs.$ We denote by $\faQp$ the positive chamber determined
by the positive system $\gS(Q)$ and by $\AQp$ the image in $\Ap$ under
the exponential map. Then $K \AQp K$ is an open dense subset of $G.$

In
view of its $\tau$-spherical behavior, the Eisenstein
integral $E(P\col \gl),$ for $P \in \minparabs,$ is
completely determined by its restriction to $\AQp.$ It follows from Harish-Chandra's result
\bib{HC2}, Thm.\ 18.1, that, on $M_0 \AQp,$ the given Eisenstein integral behaves asymptotically as follows:
\begin{equation}
\naam{e: asymptotics Eis}
E(P\col \gl \col m a)\psi \sim
 \sum_{s \in W} a^{s\gl -\rho_Q} [C_{Q|P}(s\col \gl) \psi](m), \qquad
(a \to \infty \;\;{\rm in}\; \AQp),
\end{equation}
for every $\psi \in \cAtwo,$ every $m \in M_0,$ and $\gl \in i\fapdreg.$ Here $W$
denotes the Weyl group of the root system $\gS$ and the coefficients
$C_{Q|P}(s \col \dotvar)$ are $\End(\cAtwo)$-valued analytic functions of $\gl \in i \fapdreg.$
The functions $C_{Q|P}(s \col \dotvar),$ for $s \in W,$ are uniquely determined by these
properties. A priori they have a meromorphic extension to an open neighborhood of
$i\fapd$ in $\fapdc.$

\begin{rem}
\naam{r: gl parameter C}
Harish-Chandra denotes the $c$-functions by lower case letters.
In view of Remark \refer{r: gl parameter Eis}, the $C$-functions introduced above
are related to Harish-Chandra's by the formula $C_{Q|P}(s \col i \gl) = c_{Q|P}(s \col \gl).$
\end{rem}

For $P \in \minparabs$ and $R \in \R$ we put
\begin{equation}\
\naam{e: defi fapd P R}
\fapd(P, R) := \{\gl \in \fapdc \mid \inp{\Re \gl}{\ga} < R, \;\; \forall \ga \in \gS(P) \}.
\end{equation}
A $\gS$-hyperplane in $\fapdc$ is a hyperplane of the form $\inp{\gl}{\ga} = c,$
with $\ga \in \gS$ and $c \in \C.$ The hyperplane is said to be real if
$c \in \R.$

We define $\Pi_\gS(\fapdc)$ to be the set of polynomial functions
that can be written as a product of a nonzero complex number
and linear factors of the form $\gl \mapsto \inp{\gl}{\ga} - c,$
with $\ga \in \gS$ and $c \in \C.$
The subset of polynomial functions which are products as above
with $c \in \R$ is denoted by $\Pi_{\gS, \R}(\fapdc).$

\begin{lemma}
\naam{l: estimates C functions in cone}
Let $P \in \minparabs.$
The endomorphism $C(1: \gl) = C_{P|P}(1: \gl)
\in \End(\cAtwo)$ is invertible for generic $\gl \in i\fapd.$
Both maps
\begin{equation}
\naam{e: C pm}
\gl \mapsto C(1: \gl)^{\pm 1}
\end{equation}
extend to
$\End(\cAtwo)$-valued meromorphic functions on $\fapdc$ that can be expressed
as products of functions of the form
$\gl \mapsto c_\ga(\inp{\gl}{\ga}),$
 for $\ga \in \gS(P),$
with $c_\ga$ a meromorphic function on $\C$ with real
singular locus. Accordingly, each of the functions (\refer{e: C pm})
has a singular locus equal to a locally finite union of real $\gS$-hyperplanes
in $\fapdc.$

Let $R \in \R.$ Only a finite number of the mentioned singular
hyperplanes intersect $- \fapd(P, R).$ There exist polynomial
functions  $q_\pm \in \Pi_{\gS, \R}(\fapd)$ such that $\gl \mapsto
q_\pm(\gl)  C(1\col - \gl)^{\pm 1}$ are regular on the
closure of the set $\fapd(P, R).$ If $q_\pm$ is any pair of
polynomials with these properties, there exist constants $n \in
\N$ and $C > 0$ such that
$$
\| q^{\pm}(\gl) C(1\col -  \gl)^{\pm 1} \|   \leq C (1 + |\gl|)^n, \qquad (\gl \in \fapd(P, R)).
$$
\end{lemma}

\proof
All assertions readily follow from the arguments in \bib{Arthur},
proof of Lemma 5.2, except possibly for the final estimate,
which at first follows for a particular choice of $q_{\pm}.$
A straightforward application of the Cauchy integral formula
then gives the result for arbitrary $q_{\pm}$ satisfying the hypotheses.
\qed

Following Harish-Chandra \bib{HC3}, \S 17, we define
the following normalized $C$-functions, for $P,Q \in \minparabs$ and $s \in W:$
\begin{equation}
\naam{e: defi ooC}
\ooC_{Q|P}(s\col \gl):= C_{Q|Q}(1\col s \gl)^{-1} C_{Q|P}(s \col \gl).
\end{equation}
The following result, due to Harish-Chandra, \bib{HC3},
will be of crucial importance
to us.

\begin{lemma}
\naam{l: rationality ooC}
{\rm (Harish-Chandra \bib{HC3})\ }
For all $P,Q \in \minparabs$ and $s \in W,$ the
$\End(\cAtwo)$-valued function $\gl \mapsto \ooC_{Q|P}(s \col \gl)$
has a rational extension to  $\fapdc.$

The endomorphism $\ooC_{Q|P}(s\col \gl)$
is invertible for generic $\gl \in \fapdc$
and $\gl \mapsto \ooC_{Q|P}(s \col \gl)^{-1}$ is a rational $\End(\cAtwo)$-valued
function.

Finally, each of the functions $\gl \mapsto \ooC_{Q|P}(s \col \gl)^{\pm 1}$
is a product of functions of the form
$\gl \mapsto c_\ga(\inp{\gl}{\ga}),$ for $\ga \in \gS,$
with $c_\ga$ a $\End(\cAtwo)$-valued rational function on $\C.$
\end{lemma}

\proof
The assertions for $\ooC_{Q|P}(s \col \gl)$
follow from \bib{HC3}, Lemma 19.2 combined with the
Corollary to Lemma 17.2 and with Lemma 17.4 of the same article.
For imaginary $\gl,$ the endomorphism $\ooC_{Q|P}(s \col \gl)$
is unitary, by \bib{HC3}, Lemma 17.3. The remaining assertions
now follow by application of Cramer's rule.
\qed

\eqsection{Relations between the Fourier transforms}
\naam{s: Relations between the Fourier transforms}
According to \bib{HC3}, Lemma 17.2,
the Eisenstein integrals are related by the following functional
equations:
\begin{equation}
\naam{e: functional equation uEis}
E(P\col \gl \col \dotvar) = E(Q\col s \gl \col \dotvar) \ooC_{Q|P}(s \col \gl),
\end{equation}
for $P,Q \in \minparabs$ and   $s \in W,$  as an identity
of meromorphic $\Ci(G)\otimes \Hom(\cAtwo,\Vtau)$-valued functions in the variable $\gl \in \fapdc.$

\begin{lemma}
\naam{l: inversion ooC}
Let $P,Q \in \minparabs.$ Then $\ooC_{Q|P}(1 \col \gl) \after \ooC_{P|Q}(1 \col \gl) = I,$
as an identity of $\End(\cAtwo)$-valued functions in the variable $\gl \in \fapdc.$
\end{lemma}

\proof
From the functional equation for the Eisenstein integral it
follows that, for $x \in G,$
$$
E(P\col \gl \col x) = E(P\col  \gl \col x) \ooC_{P|Q}(1 \col \gl) \ooC_{Q|P}(1 \col \gl).
$$
Using (\refer{e: asymptotics Eis}), we infer that this identity is valid with
$C_{P|P}(1 \col \gl)$ in place of $E(P\col \gl \col x)$ on both sides.
As $C_{P|P}(1 \col \gl)$ is invertible
for generic $\gl,$ the required identity follows.
\qed

In view of (\refer{e: defi unnormalized dual Eis}) it follows immediately from
(\refer{e: functional equation uEis}) that the
unnormalized dual Eisenstein integrals satisfy the following functional equations,
for $P,Q \in \minparabs$ and $s \in W,$
\begin{equation}
\naam{e: functional equation udEis}
\ooC_{Q|P}(s \col -\bar \gl)^* \;\udE(Q \col s \gl\col \dotvar ) = \udE(P\col \gl \col \dotvar),
\end{equation}
as an identity of meromorphic $\Ci(G) \otimes \Hom(\Vtau, \cAtwo)$-valued functions of $\gl \in \fapdc.$
In view of the definition of the Fourier transform in (\refer{e: defi uFou}),
this in turn implies that, for every
$f\in \Cci(G,\tau),$
\begin{equation}
\naam{e: functional equation uFou} \ooC_{Q|P}(s \col -\bar
\gl)^*\; \uFou_Q  f(s \gl) = \uFou_P f (\gl),
\end{equation}
as an identity of meromorphic $\cAtwo$-valued functions of $\gl \in \fapdc.$

\begin{lemma}
\naam{l: uAC imply transformation formula gf by C}
 Let $\gf = (\gf_P)_{P \in \minparabs} \subset
\cO(\fapdc) \otimes \cAtwo,$
and assume that $\Lau \gf = 0$ for all $\Lau \in \uAC(G,\tau, \minparabs).$
Then for all $P,Q \in \minparabs$ and $s \in W,$
\begin{equation}
\naam{e: transformation formula gf by C}
\ooC_{Q|P}(s \col -\bar \gl)^*\; \gf_Q (s \gl) = \gf_P(\gl),
\end{equation}
for generic $\gl \in \fapdc.$
\end{lemma}

\proof
Let $P,Q\in \minparabs$ and $s \in W$ be fixed.
In view of
Lemma \refer{l: rationality ooC} there exists a polynomial
function $q \in \Pi_{\gS}(\fapdc)$ such that
$\gl \mapsto q(\gl) \ooC_{Q|P}(s \col -\bar \gl)^*$
is polynomial.  Let $\gf$ fulfill the hypothesis and let $\mu \in
\fapdc\setminus q^{-1}(0).$

We define Taylor functionals
$\cL_R \in \cO(\fapdc)^*_\tayl
\otimes \cAtwo^*$ by
$\Lau_P \psi:=  - q(\mu) \psi(\mu)$ and
$\Lau_Q \psi = \ev_\mu[\gl \mapsto
q(\gl) \ooC_{Q|P}(s \col -\bar \gl)^* \psi(s \gl)],$ for $\psi \in \cO(\fapdc)
\otimes \cAtwo,$
and by $\Lau_R = 0$ for $R \in \minparabs\setminus \{P, Q\}.$
It follows from (\refer{e: functional equation udEis}) that
$\Lau \in \uAC(G,\tau, \minparabs).$ Hence, $\Lau \gf = 0.$
We conclude that $\gf$ satisfies
(\refer{e: transformation formula gf by C})
for
$\gl \in \fapdc\setminus q^{-1}(0).$
\qed

If $V$ is a finite dimensional real linear space,
we denote by $\Mer(V_\iC)$ the space of meromorphic functions
on $V_\iC.$
Given $P,Q \in \minparabs,$ we define the endomorphism
$\mult_{P|Q}$ of $\Mer(\fapdc) \otimes \cAtwo$ by
$$
[\mult_{P|Q} \psi](\gl) = \ooC_{Q| P}(1 \col - \bar \gl)^* \psi(\gl).
$$
In particular, $\mult_{P|P} = I.$

\begin{lemma}
Let $P, Q \in \minparabs$ and $R > 0.$ Then  $\mult_{P|Q}$ maps $\uPW_R(G, \tau, Q)$
continuously into $\EPW_R(\fapdc) \otimes \cAtwo.$
\end{lemma}

\proof
It follows from Lemma \refer{l: rationality ooC} that there exists
a polynomial $q \in \Pi_\gS(\fapdc)$ such that the function
$\gl \mapsto q(\gl)\, \ooC_{Q|P}(1\col -\bar \gl)^*$ is polynomial on $\fapdc.$
This in turn implies that the map $q \after \mult_{P|Q}$ maps $\uPW_R(G,\tau,Q)$
continuously into $\EPR(\fapdc) \otimes \cAtwo.$

Let $\ga \in \gS$ and $c \in \C$ be such that
$l : \gl \mapsto \inp{\gl}{\ga} - c$ is a factor of $q.$ Let $d$ be the highest
integer such that $l^d$ is still a factor of $q.$ Let $H_\ga $ denote the element
of $\fap$ determined by $H_\ga \perp \ker \ga$ and $\ga(H_\ga) = 2.$
Fix $0 \leq k < d.$ Then the element $H_\ga^k \in S(\fap),$ viewed as a constant
coefficient differential operator on $\fapdc$ satisfies $H_\ga^k q = 0$
on $l^{-1}(0).$ Fix $\gl_0 \in l^{-1}(0)$ and consider
the  functional $\cL \in \cO(\fapdc)^*_\tayl \otimes \cA_2^*$ defined by
\begin{equation}
\naam{e: Lau and j P Q}
\cL (\gf):=  \ev_{\gl_0} \after H_\ga^k [q \mult_{P|Q} \gf],
\end{equation}
for $\gf \in \cO(\fapdc) \otimes \cAtwo.$ It follows from
(\refer{e: functional equation udEis}) that, for all $x
\in G$ and $v \in \Vtau,$ the function $\mult_{P|Q} \udE(Q\col
\dotvar  \col x)v$ equals $\udE(P\col \dotvar \col
x)v,$ hence is holomorphic on $\fapdc.$ By application
of the Leibniz rule we now see that $\cL(\udE(Q\col \dotvar \col
x)v) = 0$ for all $x \in G$ and $v \in \Vtau.$ Hence,
$\cL$ belongs to $\uAChol(G,\tau, Q).$ Let now $\gf \in \uPW_R(G,
\tau, Q).$ Then it follows that (\refer{e: Lau and j P Q}) equals
zero. As this is valid for every $\gl_0 \in l^{-1}(0)$ and all $0
\leq k < d,$ it follows that $l^d$ divides $q(\gl) \mult_{P|Q}
\gf.$ Treating all factors of $q$ in this fashion, we see that
$\mult_{P|Q} \gf$ is holomorphic on $\fapdc$ outside a subset of
complex codimension $2.$ It follows that $\mult_{P|Q}$ maps
$\uPW_R(G,\tau, Q)$ into $\cO(\fapdc) \otimes \cAtwo.$ Since $q
\mult_{P|Q}$ maps $\uPW_R(G,\tau, Q)$ continuous linearly into
$\EPR(\fapdc) \otimes \cAtwo,$ it follows by a repeated
application of Cauchy's integral formula, treating the linear
factors of $q$ one at a time, that $\mult_{P|Q}$ is a continuous
linear map $\uPW_R(G,\tau, Q) \to \EPR(\fapdc) \otimes \cAtwo.$
\qed

If
 $\Lau \in \uAC(G,\tau, P),$  then the family of Taylor
functionals $(\Lau_Q')_{Q \in \minparabs}$
defined by $\Lau_P' = \Lau$ and by $\Lau_Q' = 0$ for $Q \neq P$
 belongs to $\uAC(G, \tau, \minparabs).$ Accordingly,
we may view $\uAC(G,\tau, P)$ as a subspace of $\uAC(G,\tau, \minparabs).$
It follows that $\uPW_R(G, \tau, \minparabs),$ for $R > 0,$
is a subspace of the direct sum of the spaces
$\uPW_R(G, \tau, P),$ for $P \in \minparabs.$ Moreover, by continuity
of Taylor functionals, this subspace is closed.

\begin{prop}
\naam{p: proj is topol iso}
Let $Q \in \minparabs.$ Then, for each $R > 0,$
the projection onto
the component with index $Q$ induces
a topological linear isomorphism
\begin{equation}
\naam{e: pr Q in PW}
\uPW_R(G,\tau, \minparabs)\; \too\; \uPW_R(G, \tau, Q).
\end{equation}
\end{prop}

\proof
Let $E$ denote the direct sum of a finite number of copies of
$\EPW_R(\fapdc) \otimes \cAtwo,$
labeled by the elements of $\minparabs.$  For each such element $P$
let $\pr_P: E \to \EPW_R(\fapdc) \otimes \cAtwo$ denote the projection
onto the component of label $P.$
We define the map $\mult_Q: \uPW_R(G, \tau, Q) \to E$ by
$$
\pr_P \after \mult_Q = \mult_{P|Q}, \qquad (\forall\; P \in \minparabs).
$$
Then $\mult_Q$ is continuous linear; we will show that
it maps into the subspace
$\uPW_R(G, \tau, \minparabs)$ of $E.$  Let
$(\Lau_P)_{P \in \minparabs}$ belong to $\uAChol(G,\tau, \minparabs).$
For each $P \in \minparabs$ we select $q_P \in \Pi_\gS(\fapdc)$ such that
$\gl \mapsto q_P(\gl) \ooC_{Q|P}(1 \col - \bar \gl)^*$ is a polynomial function. By Lemma
\refer{l: division Tayl by q}  below,
 there exists a
$\Lau_P'\in \cO(\fapdc)^*_\tayl \otimes \cAtwo^*$
such that
$\Lau_P = \Lau_P'\after q_P$  on $\cO(\fapdc) \otimes \cAtwo.$
By application of the Leibniz rule,
we see that
$\Lau'' = \sum_{P} \Lau'_P \after q_P \after \mult_{P|Q}$
defines an element of $\cO(\fapdc)^*_\tayl \otimes \cAtwo^*.$
 It follows from the functional equations
for the Eisenstein integral that, for every $x \in G$ and  $v_Q \in \Vtau,$
\begin{eqnarray*}
\Lau'' [\udE(Q \col \dotvar \col x)v_Q] &=& \sum_P \Lau_P'\after q_P
[ \udE(P \col \dotvar \col x)v_Q]\\
&=& \sum_P \Lau_P [\udE(P \col \dotvar \col x)v_Q] = 0.
\end{eqnarray*}
Hence, $\Lau'' \in \uAChol(G,\tau, Q).$ It follows that for $\gf\in \uPW_R(G,\tau, Q)$ we have
$$
0 = \Lau'' (\gf)= \sum_P \Lau_P' [q_P \mult_{P|Q}(\gf)].
$$
Moreover, since $\mult_{P|Q}(\gf)$  is holomorphic for each $P \in \minparabs,$ it follows that the
latter expression equals $\sum_P \Lau_P \mult_Q(\gf)_P.$ Hence, $\mult_Q$ maps $\uPW_R(G, \tau, Q)$
into $\uPW_R(G, \tau, \minparabs).$ Moreover, it does so continuously, as the latter space
carries the relative topology from $E.$

From the definition of $\mult_Q$ we see
that $\pr_Q \after \mult_Q = \mult_{Q|Q} = I.$
Moreover, if $\gf \in \uPW_R(G,\tau, \minparabs),$ then
by Lemma \refer{l: uAC imply transformation formula gf by C}  with
$s = 1,$
it follows that $\gf_P = \mult_{P|Q}(\gf_Q),$ for all $P \in \minparabs.$
Hence, $\mult_Q \after \pr_Q = I$ on
 $\uPW_R(G,\tau, \minparabs).$
It follows that $\pr_Q$ restricts to a topological linear isomorphism
(\refer{e: pr Q in PW}) with inverse $\mult_Q.$
\qed
\begin{lemma}
\naam{l: division Tayl by q}
Let $q \in \Pi_\gS(\fapdc).$ Then for every $\cL \in \cO(\fapdc)^*_\tayl$ there exists
a $\cL'\in \cO(\fapdc)^*_\tayl$ such that $\cL = \cL'\after q$ on $\cO(\fapdc).$
\end{lemma}

\proof
This may be proved in the same fashion as \bib{BSpw}, Lemma 10.5.
\qed

It is an immediate consequence of Proposition
\refer{p: proj is topol iso} that
Theorem \refer{t: arthurs PW one} and Theorem \refer{t: arthurs PW two} are equivalent.

\eqsection{The normalized Fourier transform}
\naam{s: the normalized Fourier transform}
The purpose of this section is to give equivalent versions
of the Paley--Wiener theorems discussed in the previous sections,
Theorem \refer{t: arthurs PW two} and Theorem \refer{t: distributional PW}.
The new versions are formulated in terms of a suitably normalized Fourier
transform, which in the final section will be shown
to coincide with the analogous Fourier
transform for the group viewed as a symmetric space. The normalized Fourier transform
is defined as in Section \refer{s: basic notation}, but with  a differently normalized
Eisenstein integral.

Let $P \in \minparabs.$
We define the normalized Eisenstein integral by
\begin{equation}
\naam{e: defi normalized Eis}
\nE(P\col \gl \col x) = E(P\col \gl \col x)C_{P|P}(1 \col \gl)^{-1},
\end{equation}
for generic $\gl \in \fapdc$ and for $x \in G.$
Then $\gl \mapsto \nE(P\col \gl)$ is a meromorphic $\Ci(G) \otimes \Hom(\cAtwo, \Vtau)$--valued
function on $\fapdc.$ It is not entire holomorphic anymore, but its singular set is
of a simple nature. Indeed, by Lemma \refer{l: estimates C functions in cone}
the singular set
is a locally finite union of hyperplanes of the
form $\inp{\gl}{\ga} = c,$ with $\ga \in \gS^+,$ $c \in \R.$
Moreover, the occurring
constants $c$ are bounded from below. It is known that the singular set
is disjoint from the imaginary space $i \fapd,$ but we  shall  not need
this here.

As before we define (normalized) dual Eisenstein integrals by
\begin{equation}
\naam{e: defi normalized dual Eis}
\dE(P\col \gl \col x) = \nE(P\col -\bar \gl \col x)^* \in \Hom(\Vtau, \cAtwo)
\end{equation}
for generic $\gl \in \fapdc$ and for $x \in G.$ In terms of these we
define the normalized Fourier transform
$\cF_P: \Ccminf(G\col \tau) \to \Mer(\fapdc) \otimes \cAtwo$ by
$$
\cF_P f(\gl) = \int_G \dE(P\col \gl \col x) f(x) \; dx.
$$
\begin{lemma}
Let $f \in C^{-\infty}_c(G, \tau).$ The unnormalized and normalized Fourier transforms
are related by
\begin{equation}
\naam{e: relation uFou and Fou}
C_{P|P}(1 \col - \bar \gl)^* \Fou_P f(\gl) =  \uFou_P f (\gl),
\end{equation}
as an identity of meromorphic functions of the variable $\gl \in \fapdc.$
\end{lemma}

\proof
Replacing $\gl$ by $- \bar \gl$ in both sides of (\refer{e: defi normalized Eis}),
then multiplying  with the $C$-function and taking conjugates,
we obtain, in view of (\refer{e: defi unnormalized dual Eis}) and
(\refer{e: defi normalized dual Eis}),
$$
C_{P|P}(1 \col - \bar \gl)^* \dE(P\col \gl \col x) = \udE(P\col \gl\col x),
$$
as meromorphic functions of $\gl \in \fapdc$ with values in $\Ci(G) \otimes  \Hom(\Vtau, \cAtwo).$
The result follows by testing with $f dx.$
\qed

The singular nature of the normalized Eisenstein integral does not
allow us to define Arthur--Campoli functionals in terms of Taylor functionals
as we did in Section \refer{s: Arthur's Paley--Wiener theorem}.
Instead we need the concept of  Laurent functional introduced in
 \bib{BSanfam}, Sect.\ 12. We briefly recall its definition.

Let $V$ be a finite-dimensional real linear space and let $X \subset V^*\setminus \{0\}$
be a finite subset. For a point $a \in V_\iC,$ we define the polynomial function
$\pi_a: V_\iC \to \C$  by
$$
\pi_a: = \prod_{\xi \in X} (\xi - \xi(a)).
$$
The ring of germs of meromorphic functions at $a$ is denoted by $\Mer(V_\iC,a).$
We define the subring
$$
\Mer(V_\iC, a, X): = \cup_{N \in \N}\;\; \pi_a^{-N} \cO_a.
$$
Let $\ev_a$ denote the linear functional on $\cO_a$ that
assigns to a germ $f \in \cO_a$ its value $f(a)$ at $a.$

An $X$-Laurent functional at $a \in V_\iC$ is a linear functional
$\Lau \in \Mer(V_\iC, a, X)^*$ such that for every $N \in \N$ there
exists a $u_N \in S(V)$ such that
\begin{equation}
\naam{e: expression Lau}
\Lau = \ev_a \after u_N \after \pi_a^N \qquad {\rm on} \quad \pi_a^{-N} \cO_a.
\end{equation}
The space of all Laurent functionals on $V_\iC,$ relative to $X,$ is defined
as the algebraic direct sum of linear spaces
\begin{equation}
\naam{e: defi space of Lau}
\Mer(V_\iC, X)^*_\laur:= \bigoplus_{a \in V_\iC} \Mer(V_\iC, X, a)^*_\laur.
\end{equation}
For $\Lau$ in the space (\refer{e: defi space of Lau}),
 the finite set of $a \in V_\iC$ for which
the component $\Lau_a$ is nonzero is called the support of $\Lau$ and
denoted by $\supp \Lau.$

According to the above definition, any $\Lau \in \Mer(V_\iC, X)^*_\laur$
may be decomposed as
$$
\Lau = \sum_{a \in \supp \Lau} \Lau_a.
$$
Let $\Mer(V_\iC, X)$ denote the space of meromorphic functions $\gf$ on $V_\iC$
with the property that the germ $\gf_a$ at any point $a \in V_\iC$ belongs
to $\Mer(V_\iC, a,  X).$ Then the natural bilinear map
$\Mer(V_\iC, X)^*_\laur \times \Mer(V_\iC, X) \to \C,$ given by
$$
(\Lau, \gf) \mapsto \Lau \gf:= \sum_{a \in \supp \Lau} \Lau_a \gf_a
$$
induces an embedding of $\Mer(V_\iC, X)^*_\laur$ onto a linear subspace of the dual space
$\Mer(V_\iC, X)^*.$ For more details concerning these definitions,
we refer the reader to \bib{BSanfam}, Sect.\ 12.

\begin{lemma}
\naam{l: multiplication of Lau by func}
Let $\Lau \in \Mer(V_\iC, X)^*_\laur$ and let $\psi \in \Mer(\Omega, X),$
for $\Omega$ an open neighborhood of $\supp \Lau.$ Then
$\Lau \after \psi$ belongs to $\Mer(V_\iC, X)^*_\laur.$
\end{lemma}

\proof
Without loss of generality we may assume that $\Lau$ is supported
by a single point $a \in V_\iC.$  First assume that
$\psi_a \in \cO_a.$ Then the result follows by a straightforward
application of the definition containing (\refer{e: expression Lau})
combined with the Leibniz rule. It remains to establish
the result for $\psi = \pi_a^{-k},$ with $k \in \N.$ In this
case the result is an immediate consequence of the mentioned
definition.
\qed

The following result relates the Laurent functionals to the Taylor functionals
defined in Section \refer{s: Arthur's Paley--Wiener theorem}.
The inclusion map $\cO(V_\iC) \subset \Mer(V_\iC, X)$ induces
a surjection $\Mer(V_\iC, X)^* \to \cO(V_\iC)^*.$ This property
of surjectivity also holds on the level of Laurent functionals.
The space $\cO(V_\iC)^*_\tayl$ may naturally be viewed
as a subspace of $\cO(V_\iC)^*.$ The natural map
$\Mer(V_\iC, X)^*_\laur \to \cO(V_\iC)^*$ maps into this
subspace.

\begin{lemma}
\naam{l: surjectivity natural map}
The natural map $\Mer(V_\iC, X)^*_\laur \to \cO(V_\iC)^*_\tayl$
is surjective.
\end{lemma}

\proof
Let $U \in \cO(V_\iC)^*_\tayl.$ Without loss of generality
we may assume that $U$ is supported by a single point
$a \in V_\iC.$ Then $U = \ev_a \after u$ for some $u \in S(V).$
By \bib{BSres}, Lemma 1.7 with $d' = 0,$ there exists
a $\Lau \in \Mer(V_\iC, a, X)^*_\laur$ determined
by a sequence $(u_n)_{n \in \N} \subset S(V),$ such that $u_0 = u.$
Thus, $\Lau$ restricts to $U$ on $\cO(V_\iC).$
\qed

Before proceeding, we formulate a result concerning division that will
be frequently used in the sequel.

\begin{lemma}
\naam{l: laurent and poles}
Let $E$ be a finite dimensional complex linear space.
Let $S$ be a linear subspace of $\Mer(V_\iC, X) \otimes E,$ let
$S^\circ$ be the annihilator of $S$ in $\Mer(V_\iC, X)^*_\laur \otimes E^*$
and let $S^{\circ\circ}$ be the  space of functions
$\gf \in \Mer(V_\iC, X) \otimes E$ such that $\Lau \gf = 0$
for all $\Lau \in S^\circ.$

Let $\psi$ be a nonzero  $\End(E)$-valued  meromorphic function
on $V_\iC$ such that both $\psi$ and $\psi^{-1}$ belong
to $\Mer(V_\iC, X) \otimes \End(E). $
Let $\Omega \subset V_\iC$ be an open subset such that
$\psi \gf $ is regular on $\Omega$
for all $\gf \in S.$
Then $\psi \gf$ is regular on $\Omega$ for all $\gf \in S^{\circ\circ}.$
\end{lemma}

\proof
We first assume that $\psi = I.$
Then every $\gf \in S$ is regular on $\Omega.$
Let now $\gf \in S^{\circ\circ}$ and consider a point $a \in \Omega.$
Then it suffices to show that the germ $\gf_a$ is regular at $a.$
As $\gf_a \in \Mer(V_\iC, a) \otimes E,$ there exists a
product $q$ of factors of the form $\xi - \xi(a),$ with $\xi \in X,$  such that
$q\gf_a$ is regular at $a.$ We fix $q$ of minimal degree. Then
$q \gf_a$ has a non-trivial value at $a.$ Hence, there exists a
linear functional $\eta \in E^*$ such that $\ev_a \after q \inp{\gf_a}{\eta}
\neq 0.$ There exists a Laurent functional
$\Lau \in \Mer(V_\iC, a)^*_\laur$ such that $\Lau = \ev_a$ on $\cO_a(V_\iC).$
Now $\Lau \after q$ is a Laurent functional, and it follows from
the above that  $\Lau_1:= [\Lau \after q] \otimes \eta$ is nonzero on $\gf_a.$
Hence $\Lau_1 \notin S^\circ.$ It follows that there exists a function
$\gf_1 \in S$ such that $q(a) \eta(\gf_1(a)) = \Lau_1 \gf_1$ is nonzero.
This implies that $q$ is nonzero at $a,$ hence constant. We conclude
that $\gf_a$ is regular at $a.$

We now turn to the case with $\psi$ general. Then multiplication
by $\psi$ induces a linear automorphism of $\Mer(V_\iC, X) \otimes E,$
whose inverse is multiplication by $\psi^{-1}.$
In view of Lemma \refer{l: multiplication of Lau by func},
the map $\psi^* : \Lau \mapsto \Lau \after \psi$
is a linear automorphism of $\Mer(V_\iC, X)^*_\laur \otimes E$
with inverse $(\psi^{-1})^* .$ Put $S_1: = \psi S.$ Then
$S_1^\circ = \psi^{*-1}(S^\circ)$ and $S_1^{\circ\circ} = \psi S^{\circ\circ}.$
By the first part of the proof it follows that all elements
of $S_1^{\circ\circ}$ are regular on $\Omega.$ The result follows.
\qed

We use the non-degenerate bilinear map $\inp{\dotvar}{\dotvar}$ on $\fapd$ to identify
this space with its real linear dual. Accordingly we view $\gS$ as a finite subset
of $\fa_0^{**}\setminus \{0\}$ and invoke the space of $\gS$-Laurent functionals
on $\fapdc$ in the following definition.

\begin{defi}
A (normalized) Arthur--Campoli functional for $(G, \tau,P)$ is a Laurent functional
$\cL \in \Mer(\fapdc, \gS)^*_\laur \otimes \cAtwo^*$ such that
$$
\cL\; [\dE (P\col \dotvar \col x)v] = 0
$$
for all $x \in G$ and $v \in \Vtau.$
The space of such functionals is denoted by $\AC(G,\tau,P).$
\end{defi}

Our next objective is to define suitable spaces of meromorphic functions
with controlled singular behavior.

A $\gS$-hyperplane in $\fapdc$ is defined to be a hyperplane of the form $l^{-1}(0),$
where $l: \gl \mapsto \inp{\gl}{\ga} - c$ with $\ga \in \gS$ and $c \in \C.$ The hyperplane is said
to be real if $c \in \R.$ A $\gS$-configuration in $\fapdc$ is a locally finite collection
of $\gS$-hyperplanes. The configuration is said to be real if all its hyperplanes are real.
Let now $\Hyp$ be a real $\gS$-configuration.
For each
$H \in \Hyp$ we fix $\ga_H \in \gS$ and $s_H \in \R$ such that $H$ equals the zero locus
of $l_H: \gl \mapsto \inp{\gl}{\ga} - s_H.$

Let $d: \Hyp \to \N$ be a map. For $\omega$ a subset of $\fapd$
whose closure intersects only finitely many
hyperplanes from $\Hyp,$ we define the polynomial function
$\pi_{\omega, d }$ on $\fapdc$ by
\begin{equation}
\naam{e: defi pi omega d}
\pi_{\omega,d} =
\prod_{H \in \Hyp\atop H \cap \,{\rm cl}\,\omega \neq \varnothing}l_H^{d(H)}.
\end{equation}
Moreover, we define  $\Mer(\fapdc, \Hyp, d)$ to be the space of
meromorphic functions $\gf \in \Mer(\fapdc)$ such that for every
bounded open subset $\omega$ of $\fapd,$ the function
$\pi_{\omega, d} \gf$ is regular on $\omega + i\fapd.$

For $\omega$ a bounded  subset of $\fapd$ and $n \in \Z$ we
define the $[0,\infty]$-valued seminorm $\nu_{\omega, d, n}$ on
$\Mer(\fapdc, \Hyp, d)$
by
\begin{equation}
\naam{e: seminorm cP space}
\nu_{\omega, d, n} (\gf)
: = \,\sup_{\gl \in \omega  +   i \fapd}
 \;( 1 + |\gl|)^n \, |\pi_{\omega, d}(\gl)\gf(\gl)|.
\end{equation}
We define
$
\cP(\fapdc, \Hyp, d)
$
to be the space of functions $\gf \in \cM(\fapdc, \Hyp, d)$ such that
$
\nu_{\omega, d, n}(\gf) < \infty
$
for every compact set $\omega \subset \fapd$ and all $n \in \N.$ This space
is a Fr\'echet space with topology induced by the collection of
seminorms $\nu_{\omega, d, n},$ for $\omega$ compact and $n \in \N.$

We denote by $\cN = \cN(\fapd)$ the collection of
maps $n: \cC \to \N,$ with $\cC$ the collection of compact subsets
of $\fapd.$ On $\cN$ we define the partial ordering
$\leq$ by $n \leq m$
if and only if $n(\omega) \leq m(\omega)$ for all $\omega \in \cC.$
For $n \in \cN$ we define
$\cP_n^*(\fapdc, \Hyp, d)$ to be the space of functions
$\gf \in \Mer(\fapdc, \Hyp, d)$
such that
$$
\nu_{\omega,d, - n(\omega)} (\gf) = 
\sup_{\lambda \in \omega + i{\mathfrak a}_0^*} (1 + |\lambda|)^{-n(\omega)} \,
|\pi_{\omega, d}(\lambda) \varphi(\lambda)|
 < \infty
$$
for every compact subset
$\omega \subset \fapd.$  Equipped with the seminorms
$\nu_{\omega, d, - n(\omega)}$ the space $\cP_n^*(\fapdc, \Hyp, d)$
is a complete locally convex space. If  $m \leq n$
then clearly $\cP_m^* \subset \cP_n^*,$
with continuous linear inclusion
map.

We now define
$$
\cP^*(\fapdc, \Hyp, d)  :=  \cup_{n \in \cN(\fapd)}\; \; \cP_{n}^* (\fapdc, \Hyp, d),
$$
and equip this space with the inductive limit locally convex topology.

In particular, these definitions may be interpreted for $\Hyp = \void$
and $d = \void.$ In this case we have $\pi_{\omega, d} = 1$
for every $\omega \subset \fapd,$ so that
$$
\cP(\fapdc, \varnothing) \quad{\rm and} \quad \cP^*(\fapdc, \void)
$$
are just spaces of holomorphic functions $\gf$ on $\fapdc$
determined in the above fashion by seminorms of the form
$$
\nu_{\omega, n}(\gf) := \sup_{\gl \in \omega + i\fapd} \;
(1 + |\gl|)^n |\gf(\gl)|.
$$
For the rest of this section,
let $P \in \minparabs$ be fixed and let $\Hyp = \Hyp_{G,\tau, P}$
be the smallest collection of $\gS$-hyperplanes in
$\fapdc$ such that the singular locus of
$\gl \mapsto \dE(P\col \gl \col \dotvar)$
is contained in the union $\cup \Hyp.$ In view of Lemma
\refer{l: estimates C functions in cone} the collection $\Hyp$ is locally finite and
consists of real $\gS$-hyperplanes.
Moreover, by the same lemma, the set
of $H \in \Hyp$ with $H \cap \fapd(P,R) \neq \void$ is finite,
for every $R \in \R.$

We define the map $d= d_{G,\tau, P}: \Hyp \to \N$
as follows. For each $H \in \Hyp$ we fix $l_H: \fapdc \to \C$
as in (\refer{e: defi pi omega d}) and define
$d(H)$ as the smallest integer $k \geq 0$ such
that the $\Ci(G) \otimes \Hom(\Vtau, \cAtwo)$-valued
meromorphic function $l_H ^k \dE(P\col \dotvar)$
extends regularly over $H \setminus \cup \{H'\in \Hyp \mid H'\neq H\}.$

Given a subset $\omega \subset \fapd$ whose closure
meets only finitely many
hyperplanes from $\Hyp,$ we define the polynomial function $\pi_{\omega, d}$
as in (\refer{e: defi pi omega d}).
In particular, we write $\pi = \pi_P$ for this
polynomial with $\omega = \fapd \cap \bfapd(P,0),$ where the second set in
the intersection denotes
the closure of
the set (\refer{e: defi fapd P R}) with $R = 0.$
Thus, the $\Ci(G) \otimes \Hom(\Vtau, \cAtwo)$-valued
function $\gl \mapsto \pi(\gl) \dE(P\col \gl)$ is holomorphic
on a neighborhood of  $\bfapd(P, R)$
and $\pi \in \Pi_{\gS, \R}(\fapd)$ is minimal with this property.

We define the following closed subspace  of
$\cP(\fapdc, \Hyp, d) \otimes \cAtwo:$
\begin{eqnarray}
\lefteqn{\cP_\AC(\fapdc, \Hyp, d, P)}\nonumber
\\
&:=&
\{ \gf \in \cP(\fapdc, \Hyp, d) \otimes \cAtwo \mid  \Lau \gf  = 0,\;
\forall \Lau \in \AC(G,\tau, P) \}
\end{eqnarray}
Finally, we define the space $\cP_\AC^*(\fapdc, \Hyp, d, P)$
in a similar fashion, but with $\cP$
replaced by $\cP^*.$
\pagebreak
\begin{defi}
\naam{d: normalized PW}
\begin{enumerate}
\itema
Let $R > 0.$ We define the {\it Paley--Wiener space} $\PW_R(G,\tau,P)$ to be the subspace of
$\cP_\AC(\fapdc, \Hyp, d,P)$ consisting of functions $\gf$ such that, for all $n \in \N,$
\begin{equation}
\naam{e: PW estimate}
\sup_{\gl \in \bar\fapd(P, 0)} (1 + \|\gl\|)^n e^{- R |\Re \gl|} \|\pi(\gl) \gf(\gl)\| < \infty.
\end{equation}
The space is equipped with the relative topology.
\itemb
For $R > 0$ we define the {\it distributional Paley--Wiener space}
$\PW^*_R(G,\tau, P)$ to be the subspace
of $\cP_\AC^*(\fapdc, \Hyp, d,P )$ consisting of functions $\gf$ for which there exists
a constant $n \in \N$ such that
\begin{equation}
\naam{e: PW star estimate}
\sup_{\gl \in \bar\fapd(P, 0)} (1 + \|\gl\|)^{-n } e^{- R |\Re \gl|} \|\pi(\gl) \gf(\gl)\| < \infty.
\end{equation}
This space is also equipped with the relative topology.
\end{enumerate}
\end{defi}

We will finish this section by  discussing the relation of these
Paley--Wiener spaces  with the unnormalized Paley--Wiener spaces
introduced in Definitions \refer{d: unnormalized PW}
and \refer{d: unnormalized distrib PW}.
As a preparation,
we first give another characterization of the unnormalized Paley--Wiener spaces.

We define
$
\cP_{\uAC}(\fapdc, \void, P)
$
and
$
\cP_{\uAC}^*(\fapdc, \void, P)
$
as the closed subspaces of  the spaces $\cP(\fapdc, \void) \otimes \cAtwo$
and $\cP^*(\fapdc, \void) \otimes \cAtwo,$ respectively, consisting
of the functions $\gf$ satisfying the relations $\Lau \gf = 0$ for
all $\Lau \in \uAChol(G, \tau, P).$

\begin{prop}
\naam{p: uPW by estimate on cone}
Let $R > 0.$
\begin{enumerate}
\itema
The space
$\uPW_R(G, \tau, P)$ consists of the functions
$\gf \in \cP_{\uAC}(\fapdc, \void, P)$
with the property that, for every $n \in \N,$
\begin{equation}
\naam{e: uPW estimate}
\sup_{\gl \in \bar\fapd(P, 0)} (1 + |\gl|)^n e^{-R |\Re \gl|}
\|\gf(\gl)\| < \infty.
\end{equation}
Moreover, the topology of $\uPW_R(G, \tau, P)$ coincides with the
relative topology from $\cP_{\uAC}(\fapdc, \void, P).$
\itemb
The space
$\uPW_R^*(G, \tau, P)$ consists of the functions
$\gf \in \cP_{\uAC}^*(\fapdc, \void, P)$
with the property that there exists a number $n \in \N$ such that
\begin{equation}
\naam{e: uPW star estimate}
\sup_{\gl \in \bar\fapd(P, 0)} (1 + |\gl|)^{-n} e^{-R |\Re \gl|}
\|\gf(\gl)\| < \infty.
\end{equation}
Moreover, the topology of $\uPW_R^*(G, \tau, P)$ coincides with the
relative topology from $\cP_{\uAC}^*(\fapdc, \void, P).$
\end{enumerate}
\end{prop}

\proof We start by making some remarks on Euclidean Paley--Wiener
spaces. From the text preceding Definition \refer{d: arthurs PW}
we recall the definition of the space $\EPW_R(\fapdc),$ equipped
with the Fr\'echet topology $\cT$ induced by the seminorms
$\nu_{R,n},$ for $n \in \N.$ Clearly, $\EPW_R(\fapdc) \subset
\cP(\fapdc,  \void),$ with continuous inclusion map. We denote by
$\cT_r$ the associated relative topology on $\EPW_R(\fapdc).$ Then
$\cT$ is finer than $\cT_r.$ We will show that both topologies are
in fact equal.

By the Euclidean Paley--Wiener
theorem, Euclidean Fourier transform $\eFou$
defines a continuous linear isomorphism from $C^\infty_R(\fap)$ onto
$\EPW_R(\fapdc).$
From a straightforward estimation it follows
that the inverse Fourier transform $\eFou^{-1}$ is continuous
from $\cP(\fapdc, \void)$ to $C^\infty(\fap).$
It follows from the above that the identity
map $\eFou \after \eFou^{-1}$
is continuous from $(\EPW_R(\fapdc),\cT_r)$ to
$(\EPW_R(\fapdc), \cT).$ Hence $\cT = \cT_r.$

From the text
preceding Definition \refer{d: unnormalized distrib PW}
we recall the definition of $\EPW^*_R(\fapdc),$
 equipped with the inductive limit locally convex topology
denoted $\cT^*.$
Clearly,
$\EPW^*_R(\fapdc) \subset \cP^*(\fapdc, \void),$
with continuous inclusion map. Let $\cT^*_r$ denote the associated relative
topology on  $\EPW^*_R(\fapdc).$
Then $\cT^*$ is finer than $\cT^*_r.$
We will show that both topologies are equal.

By the distributional Euclidean Paley-Wiener
theorem, $\eFou$ maps $C^{-\infty}_R(\fap)$
bijectively onto  $\EPW_R^*(\fapdc).$
Fix $R' > R$ and let
$k$ be an arbitrary positive integer.
Then  $C^{-\infty}_R(\fap)_k,$ the subspace of generalized
functions of order at most $k,$ naturally embeds into
the continuous linear dual of the Banach space $C^k_{R'}(\fap),$
equipped
with the $C^k$ norm $\|\,\cdot\,\|_{C^k}.$ Accordingly,
we equip $C^{-\infty}_R(\fap)_k$
with the restriction of the dual norm.
By a straightforward estimation,
there exists a $C_k > 0$ such that
$$
\nu_{R, k}(\eFou(f)) \leq C_k \|f\|_{C^k}
$$
for all $f \in C^k_{R'}(\fap).$ Let $n \in \cN(\fapd).$ Then by transposition
it readily follows that
$\eFou^{-1}$   maps
$\cP_{n}^*(\fapdc, \void) \cap \EPW_R^*(\fapdc)$
into $C^{-\infty}_R(\fap)_{k},$
 with $k = n(\{0\}) +  \dim \fap + 1.$
Moreover, this map is continuous with respect
to the relative topology from $\cP_{n}^*(\fapdc, \void)$
on the first of these
spaces. Now $\eFou$ maps
$C^{-\infty}_R(\fap)_{k}$ continuously into
$\EPW_R^*(\fapdc).$
It follows that the identity map $\eFou\after
\eFou^{-1}$
is continuous from
$\cP^*_n(\fapdc,\void) \cap \EPW_R^*(\fapdc)$
to $\EPW_R^*(\fapdc).$
By the universal property
of the inductive limit, it follows
that $\cT_r^*$ is finer than $\cT^*.$ Hence $\cT^* = \cT_r^*.$

We proceed with the actual proof.
We denote the subspace of $\cP_{\uAC}(\fapdc, \void, P)$ defined
in (a) by $\cP\cW$ and the similar subspace defined in (b) by $\cP\cW^*.$
These subspaces are equipped with the relative topologies.
Clearly, $\uPW_R(G, \tau, P)$ is a subspace of
$\cP\cW,$ and $\uPW^*_R(G,\tau, P)$ a subspace of $\cP\cW^*,$
with continuous inclusion maps.
To conclude the proof we must
establish the converse inclusions, also with continuous inclusion maps.

A holomorphic function
$\gf \in \cO(\fapdc) \otimes \cAtwo$ that is annihilated by
$\uAChol(G, \tau, P),$ satisfies the functional
equations
\begin{equation}
\naam{e: functional equation phi}
\gf(\gl) =  \ooC_{P|P}(s \col -\bar \gl)^*\gf(s\gl),
\end{equation}
for $s \in W$ and generic $\gl \in \fapdc;$ in view of Proposition
\refer{p: proj is topol iso} this follows from Lemma \refer{l: uAC imply transformation formula gf by C}
with $P = Q.$

In view of Lemma \refer{l: rationality ooC}, there exists a product $q_s$ of linear factors
of the form $\inp{\dotvar}{\ga} - c,$ with $\ga \in \gS$ and $c \in \C,$
such that $\gl \mapsto q_s(\gl) \ooC_{P|P}(s \col -s^{-1} \bar \gl)^*$ is polynomial.
We define
$$
q(\gl) : = \prod_{s \in W} q_s(s\gl), \qquad (\gl \in \fapdc).
$$
Then there exist constants $C > 0$ and $N \in \N,$  such
that, for all $s \in W$  and $\gl \in \fapdc,$
$$
\|q(s^{-1}\gl) \ooC_{P|P}(s \col - s^{-1} \bar \gl) \| \leq C (1 + |\gl|)^N.
$$
If we combine this with the functional equation (\refer{e: functional equation phi}),
we see that,
for each $s \in W,$ every $n \in \Z$ and all $\gl \in \bar\fapd(P, 0),$
$$
( 1+ |s^{-1} \gl|)^n e^{- R |\Re s^{-1}\gl|} \|q(s^{-1} \gl)\gf(s^{-1}\gl)\|
 \leq C ( 1 + |\gl|)^{n + N} e^{- R |\Re \gl|} \|\gf(\gl)\|.
$$
Combining these estimates for $s \in W,$ we obtain
$$
\nu_{R, n}(q\gf)
\leq C \sup_{\gl \in \bar\fapd(P,0)}
( 1 + |\gl|)^{n + N} e^{- R |\Re \gl|} \|\gf(\gl)\|,
$$
where $\nu_{R,n}$ is defined as in the first part of the proof.
On the other hand, by an easy application of Cauchy's integral
formula it follows that for every $n \in \Z$ there exists a
constant $C_n > 0$ such that
$$
\nu_{R,n}(\gf) \leq C_n \nu_{R,n}(q \gf),
$$
for all $\gf \in \cO(\fapdc) \otimes \cAtwo.$
It follows from these estimates
that $\cP\cW$ equals the intersection of
$\cP_{\uAC}(\fapdc, \void, P)$
with
the Euclidean Paley--Wiener space
$\EPW_R(\fapdc)\otimes \cAtwo.$
By definition, the topology of $\cP\cW$ equals
the relative topology from the first of
these spaces.
By the first part of the proof, the topology also coincides
with the relative topology from the second
of these spaces.
It follows that $\cP\cW \subset \uPW_R(G,\tau, P)$
with continuous inclusion map.
This establishes (a). Assertion (b) follows by a similar argument.
\qed

We define the map $ \multCP \in \Aut(\Mer(\fapdc) \otimes \cAtwo)$
by
$$
\multCP \gf (\gl) = C_{P|P}(1 \col - \bar \gl)^* \gf(\gl).
$$
Then  (\refer{e: relation uFou and Fou}) may be rephrased as
$$
\uFou_P= \multCP \after \Fou_P.
$$
\begin{thm}
\naam{t: multC is iso PW spaces}
Let $R >0.$ The map $\multCP \in \Aut(\Mer(\fapdc) \otimes \cAtwo)$
restricts
to a topological linear isomorphism
$$
\PW_R^*(G,\tau,P) \;\; \pijl{\simeq} \;\; \uPW^*_R(G,\tau,P),
$$
and similarly to a topological linear isomorphism
$$
\PW_R(G,\tau,P)\;\; \pijl{\simeq}\;\; \uPW_R(G,\tau,P).
$$
\end{thm}

\proof
It follows from Lemma \refer{l: rationality ooC}
that the functions $\gl \mapsto C_{P|P}(1: - \bar \gl)^{*\pm 1}$
belong to $\Mer(\fapdc, \gS) \otimes \End(\cAtwo),$ so that
the map $\multCP$ restricts to a linear automorphism
of the space $\Mer(\fapdc, \gS)\otimes \cAtwo.$
It follows from
Lemma \refer{l: multiplication of Lau by func} that transposition
induces an automorphism $\multCPt$ of
$\Mer(\fapdc, \gS)^*_\laur \otimes \cAtwo^*.$

Let $\uAC(G,\tau,P)$ denote the space of Laurent functionals
$\Lau \in \Mer(\fapdc, \gS)^*_\laur \otimes \cAtwo^*$
such that (\refer{e: Lau on udE}) holds for all $x \in G$ and $v \in \Vtau.$
Then it follows from Lemma \refer{l: surjectivity natural map} that
the natural map
\begin{equation}
\naam{e: natural map to uAChol}
\uAC(G,\tau,P) \to \uAChol(G,\tau,P)
\end{equation}  is surjective.

If $\Lau \in \Mer(\fapdc, \gS)^*_\laur\otimes \cAtwo^*,$ then
$$
\Lau [\udE(P\col \dotvar \col x)v] =
\Lau \after \multCP[ \dE(P\col \dotvar \col x)v],
$$
for all $x \in G$ and $v \in \Vtau.$
It follows that $\multCPt$ restricts to a linear isomorphism
from $\uAC(G,\tau,P )$ onto the space $\AC(G, \tau,P ).$

Let $\Mer_\AC(\fapdc, \gS, P)$ denote the space of $\gf \in \Mer(\fapdc, \gS) \otimes \cAtwo$
such that $\Lau \gf = 0$ for all $\Lau \in \AC(G,\tau, P).$
Similarly, let
$\Mer_{\uAC}(\fapdc, \gS, P)$ denote the space of $\gf \in \Mer(\fapdc, \gS) \otimes \cAtwo$
such that $\Lau \gf = 0$ for all $\Lau \in \uAC(G,\tau, P).$
Then it follows
from the above that $\multCP$ defines a linear isomorphism from
$\Mer_\AC(\fapdc, \gS, P)$ onto $\Mer_{\uAC}(\fapdc, \gS,P).$
The latter of the two spaces consists of holomorphic functions,
by  Lemma \refer{l: laurent and poles},
applied with $S$ consisting of the functions $ \udE(P\col \dotvar \col x)v,$ for
 $x \in G$ and $v \in \Vtau.$
By surjectivity of the map (\refer{e: natural map to uAChol}) it follows that
$\Mer_{\uAC}(\fapdc, \gS, P)$ equals the space $\cO_{\uAC}(\fapdc, P)$
 of $\gf \in \cO(\fapdc) \otimes \cAtwo$ with
$\Lau \gf = 0$ for all $\Lau \in \uAChol(G,\tau, P).$

We conclude that
$\multCP$ defines a linear isomorphism from $\Mer_\AC(\fapdc, \gS,P)$
 onto
$\cO_{\uAC}(\fapdc, P).$ Only the estimates remain to be taken care of.

Let $\omega \subset \fapd$ be a bounded open subset.
Applying Lemma \refer{l: laurent and poles} with $S$ consisting of all functions
$\dE(P\col \dotvar \col x),$ for $x \in G$ and with
$\Omega =  \omega  +  i \fapd$ and
$\psi = \pi_{\omega, d} \otimes I,$
we see that for every $\gf \in \Mer_\AC(\fapdc, \gS,P )$
the function $\pi_{\omega, d} \gf$ is holomorphic on
$ \omega  +  i\fapd.$ This implies
that $\Mer_\AC(\fapdc, \gS, P) \subset \Mer(\fapdc, \Hyp, d) \otimes \cAtwo.$

Again, let $\omega \subset \fapd$
be a bounded open subset. Then
$ \omega  +   i\fapd \subset \fapd(P,r)$ for a suitable real number $r.$
By Lemma \refer{l: estimates C functions in cone}, there exists a polynomial $p \in \Pi_\gS(\fapd)$
such that $\gl \mapsto p(\gl) C_{P|P}(1\col - \bar \gl)^*$ is holomorphic on $\fapd(P,r).$
Moreover, by application of the same lemma, there exist $N \in \N$ and $C > 0 $
such that, for every $n \in \Z,$
$$
\nu_{\omega,  n}(p\,\pi_{\omega, d} \multCP \gf)
\leq C \nu_{\omega, n +  N}(\pi_{\omega, d}\gf)
$$
for all $\gf \in \Mer(\fapdc, \Hyp, d) \otimes \cAtwo.$
On the other hand, let $\omega_0$ be a relatively compact subset of $\omega.$ Then,
by an easy application of Cauchy's integral formula, there exists
for every
$n \in \Z$ a constant $C_n >0,$ such that
$$
\nu_{\omega_0, n}(\psi) \leq
C_n \nu_{\omega, n}(p\pi_{\omega, d}\psi)
$$
for all $\psi \in \Mer(\fapdc, \Hyp, d) \otimes \cAtwo$ that are regular on
 $\omega  +  i\fapd.$
It follows that
$$
\nu_{\omega_0, n}(\multCP \gf)  \leq C_n C\, \nu_{\omega,d, n+ N} (\gf),
$$
for all $\gf \in \cP_\AC(\fapdc, \Hyp, d, P).$
This implies that $\multCP$ maps
$\cP_\AC(\fapdc, \Hyp, d,P )$
continuous linearly into $\cP(\fapdc, \void) \otimes \cAtwo,$
hence also continuously into $\cP_{\uAC}(\fapdc, \void, P).$ Moreover, the
same statement holds for the spaces with superscript $*.$
By a similar argument,
involving Lemma \refer{l: estimates C functions in cone}
for the inverse of the $C$-function,
it follows that $\multCP^{-1}$ maps $\cP_{\uAC}(\fapdc, \void, P)$
continuous linearly into $\cP_\AC(\fapdc, \Hyp, d, P).$
A similar statement is true
for the spaces with the superscript $*.$

Finally, using Lemma \refer{l: estimates C functions in cone} once more in the above fashion,
it follows that a function
$\gf \in \cP_\AC(\fapdc, \Hyp, d,P)$
satisfies the estimate (\refer{e: PW estimate})
for all $n \in \N$ if and only if $\multCP \gf$ satisfies the estimate
(\refer{e: uPW estimate}) for all  $n \in \N.$ If we combine this
with Proposition \refer{p: uPW by estimate on cone} (a),
we see that $\multCP$ restricts
to a topological linear isomorphism from $\PW_R(G,\tau, P)$ onto
$\uPW_R(G, \tau, P).$
The analogous  statements for the spaces with the superscript $*$ are proved
in a similar fashion.
\qed

\eqsection{The group as a symmetric space}
We retain the notation of the previous sections.
In this section we will view the group $G$ as a symmetric space, and compare the
Fourier transforms and Paley--Wiener spaces for $G$ with those for the
associated symmetric space. This will allow us to deduce the Paley--Wiener
theorems for the group from the analogous theorems for symmetric spaces.

As $G$ is of the Harish-Chandra class, the group $\spG:= G \times G$ is of this class as well.
We consider the involution
$\spgs$ of $\spG$ defined by $\spgs(x,y) = (y,x).$ Its group of fixed points, the diagonal subgroup,
is denoted by $\spH.$ The space $\spspX:= \spG / \spH$ is a reductive symmetric space of
the Harish-Chandra class.
The map $\spG \to G$ given by $(x,y) \mapsto xy^{-1}$ induces a
diffeomorphism
\begin{equation}
\naam{e: defi p on spaces}
p: \;\;\spspX:= \spG/\spH \;\longrightarrow \;G,
\end{equation}
 intertwining the natural left action of $\spG$ with the
action of $\spG$ on $G$ given by $(x,y) g = x g y^{-1}.$ Accordingly,  $G$ becomes a reductive
symmetric space of the Harish-Chandra class. We fix a choice of Haar measure $dg$ on $G;$
then $dx = p^*(dg)$ is a choice of $\spG$-invariant measure on $\spspX.$

The map $\spCartan := (\Cartan, \Cartan)$ is a Cartan involution
of $\spG$ which commutes with $\spgs.$ The associated maximal compact subgroup equals $\spK := K \times K.$
We recall that $\tau = (\tau_1, \tau_2)$ is a double unitary representation of $K$ in $\Vtau$
and define the unitary representation of $\spK$ in $\Vtau$ by $\sptau(k_1, k_2)v =\tau(k_1)v \tau(k_2)^{-1}.$
Then pull-back by $p$ induces a topological linear  isomorphism
\begin{equation}
\naam{e: sp iso for functions on G}
p^*:\;\; C^\minf(G\col \tau) \;\; \pijl{\simeq} \;\;
C^\minf(\spspX\col \sptau) ,
\end{equation}
which we shall also denote by $f \mapsto \spf.$
Clearly this isomorphism restricts to an isomorphism between the subspaces indicated
by $\Ccminf, \Ci$ and $\Cci.$

The $-1$ eigenspaces of $\spCartan$ and $\spgs$ in $\spfg$ equal $\spfp: = \fp \times \fp$
and  $\spfq:= \{(X, -X) \mid X \in \fg\},$ respectively. It follows that a maximal abelian
subspace $\spfaq$ of $\spfp \cap \spfq$ is given by
$$
\spfaq:= \{(X, - X) \mid X \in \fap \}.
$$
The derivative of $p$ equals the isomorphism $\spfg/\spfh \to \fg$ induced
by the map $\fg \times \fg \to \fg, \; (X,Y) \mapsto X - Y;$
we will denote this derivative
by $p$ as well. The map $p$ restricts to the isomorphism from $\spfaq$ onto $\fap$
given by $(X, -X) \mapsto 2X.$ Via
pull-back under the
isomorphism $p,$ we transfer the given inner product on $\fap$
to an inner product on $\spfaq.$ Accordingly, for every
$R>0$ the closed ball $\spbarB_R$ of center $0$ and radius $R$
in $\spfaq$ is mapped onto the similar ball
$\barB_R $ in $\fap.$  It follows that $p: \spspX \to G$
maps
$\spspX_R:= \spK \exp \spbarB_R \spH$ onto $G_R= K \barB_R K.$
The following result is now obvious.

\begin{lemma}
\naam{l: sp iso on R supported functions}
Let $R > 0.$
The map $p^*,$ defined in (\refer{e: sp iso for functions on G}), restricts to
a topological linear isomorphism from
$C^{- \infty}_R(G\col \tau)$ onto $C^{- \infty}_R(\spspX\col \sptau),$
and, similarly, to a topological linear isomorphism from
$C^{\infty}_R(G\col \tau)$ onto $C^{\infty}_R(\spspX\col \sptau).$
\end{lemma}

We will now compare the definition of the normalized Eisenstein
integral for $(\spspX, \sptau)$ given in \bib{BSft}, Sect.\ 2, with the one for
$(G,\tau)$ given in the present paper.
The isometry $p: \spfaq \to \fap$ induces an isometry
$p^*: \fapd \to \spfaqd$ (for the dual inner products on these
spaces). The complex linear extension of this map
is denoted by $p^*: \gl \mapsto \spgl,$ $\fapdc \to \spfaqdc.$

The system $\spgS$ of restricted roots of $\spfaq$ in $\spfg$ consists
of the roots $\frac12 \spga,$ for $\ga \in \gS.$ The root space for the root
$\frac12 \spga$ is given by $\fg_\ga \times \{0\} \oplus \{0\} \times \fg_{-\ga}.$
Thus, if $P\in \minparabs$ then $\spP:= P \times \bar P$ belongs
to the set $\spPmin$ of minimal $\spgs\spCartan$-stable parabolic subgroups of $\spG$ containing $\spAq;$
the associated system of positive roots is $\spgS(\spP) = \{ \frac 12 \spga \mid \ga \in \gS(P)\}.$
As usual, let
$
\rho_{\spP} \in \spfaqd$ be defined by
$$
\rho_{\spP}(\dotvar): = \frac 12 \tr (\ad(\dotvar)|\fn_{\spP}).
$$
Then $\rho_{\spP} = \spaceind(\rho_P).$ Thus, without ambiguity, we may use the notation
$\sprho_P$ for this functional.

Via the isometry  $\spfaq \simeq \fap$ we see that every element of the Weyl
group of $\spgS$ can be realized by an element of $\spK \cap \spH = \diag(K).$ It follows
that the coset space $\spW/\spW_{\spK \cap\spH}$ consists of one element. Thus as a set of representatives
for this coset space in the normalizer of $\spfaq$ in $\spK$ we may fix $\spcW = \{e\}.$
Accordingly, the space $\oC(\sptau)$ of \bib{BSft}, Eqn.~(17), now denoted by $\spcAtwo,$
is given by
$$
\spcAtwo:  =  \Ci(\spM/\spM\cap \spH \col \sptau) = L^2(\spM/\spM\cap \spH \col \sptau).
$$
We equip the space $\spM/\spM \cap \spH$ with the pull-back of the invariant measure
on $M_0$ under the analogue of the map (\refer{e: defi p on spaces})
for the tuple $(M_0, \tau_0),$ and
the space $\spcAtwo$ with the associated $L^2$-type inner product.
Then the  analogue of the isomorphism (\refer{e: sp iso for functions on G})
 for the tuple
$(M_0 , \tau_0)$ gives a unitary isomorphism
$$
p^*:\;\psi \mapsto \pspsi,\;\; \cA_2 \to \spcAtwo.
$$
Given $\psi \in \cA_2,$ we define the Eisenstein integral $E(\spP\col \sppsi \col \spgl)$
as in \bib{BSft}, Eqn.\ (20).
Then we have the following relation with the Eisenstein integral defined in (\refer{e: defi unnormalized E}).

\begin{lemma}
\naam{l: comparison Eis}
Let $P \in \minparabs,$ $\psi \in \cAtwo.$ Then for every $(x,y) \in \spG,$
\begin{equation}
\naam{e: comparison Eis}
E(\spP \col \sppsi \col \spgl )(x, y) =  E(P \col C_{\bar P|P} (1: - \bar \gl)^* \psi \col \gl)(xy^{-1}),
\end{equation}
as an identity of meromorphic functions in the variable $\gl \in \fapdc.$
\end{lemma}

\proof
We briefly write $N = N_P.$
Let $\gl \in \fapdc$ be such that $\Re \gl + \rho_P$ is $\bar P$-dominant. Then  $\Re \spgl + \sprho_P$
is $\spbarP$-dominant. Let $\sptpsi(\spgl): \spG \to V_\tau$ be defined as in
\bib{BSft}, Eqn.\ (17), for the situation at hand. Then $\sptpsi(\spgl) = 0$ outside $\spP \spH$ and
$$
\sptpsi(\spgl \col n a m_1 g\,, \,\bar n a^{-1} m_2 g)  = a^{2 \gl + 2 \rho_P} \psi(m_1 m_2^{-1}),
$$
for $n \in N, \bar n \in \bar N, a \in \Ap, m_1,m_2 \in M_0$ and $g \in G.$ It follows that
$\sptpsi(\spgl) = \spaceind[\tpsi(\gl)],$ where $\tpsi(\gl): G \to \Vtau$ is defined
to be zero outside $N \Ap M_0 \bar N$ and
$$
\tpsi(\gl \col nam \bar n) = a^{\gl + \rho_P} \psi(m)
$$
for $(n,a,m,\bar n ) \in N \times \Ap \times M_0 \times \bar N.$
 In view of \bib{BSft}, Eqn.\ (20),
we now infer that
\begin{eqnarray}
\lefteqn{ E(\spP\col \sppsi\col  \spgl \col (x_1, x_2)) =} \nonumber\\
 &= &
\int_{K \times K} \sptau(k_1, k_2)^{-1} \sptpsi(\spgl \col k_1 x_1 , k_2 x_2)
                      \; dk_1\, dk_2  \nonumber \\
&=&
\int_{K \times K} \tau(k_1)^{-1} \tpsi(\gl \col k_1 x_1 x_2^{-1} k_2^{-1}) \tau(k_2)
                      \; dk_1\, dk_2\nonumber \\
\naam{e: Eis with Psi}
&=&
E(P \col \Psi(\gl) \col \gl \col x_1 x_2^{-1}),
\end{eqnarray}
see (\refer{e: defi unnormalized E}),
where the function $\Psi(\gl): \;M_0 \to \Vtau$ is defined by
\begin{equation}
\naam{e: int for Psi}
\Psi(\gl \col m) = \int_K \tpsi(\gl \col m k^{-1})\tau(k) \; dk.
\end{equation}
As $M_0$ normalizes $\bar N,$ the function $\tpsi(\gl)$ transforms
according to $\tpsi(\gl\col xm) = \tpsi(\gl \col x) \tau(m),$
for $x \in G$ and $m \in M_0.$ It follows that the integrand in (\refer{e: int for Psi})
is a left $M_0$-invariant measurable function on $K.$ We now consider the real analytic
map $(\kappa, H, \nu): G \to K \times \fap \times N$ determined by
\begin{equation}
\naam{e: Iwasawa deco of element}
x = \kappa(x) \exp H(x) \nu(x),\qquad (x \in G).
\end{equation}
Then the Haar measure $d \bar n$ may be normalized such that for every $\gf \in C(K/M_0 )$ we have
$$
\int_K \gf(k) \; dk = \int_{\bar N} \gf(\kappa(\bar n))\, e^{- 2 \rho_P H(\bar n)}\; d \bar n.
$$
We apply this substitution of variables to the integral (\refer{e: int for Psi}).
Since  $\kappa(\bar n) = \bar n \nu(\bar n)^{-1} \exp [- H(\bar n)],$
whereas $\tpsi(\gl)$ is right $\bar N$-invariant, it follows that
\begin{equation}
\naam{e: Psi and C}
\Psi(\gl, m) = \int_{\bar N} e^{\inp{\gl - \rho_P}{ H(\bar n)}}\psi(m) \tau(\kappa(\bar n)) \; d\bar n
= [ C_{\bar P| P}(1 \col -\bar \gl)^* \psi](m).
\end{equation}
The last equality follows from \bib{HC2}, \S 19, Thm.\ 1, since $\tau$ is unitary
(take Remark \refer{r: gl parameter C} into account).
In the notation of \bib{HC2}, we have
$\mu_P = 0$ since  $P$ is minimal.

Combining (\refer{e: Eis with Psi}) with (\refer{e: Psi and C}) we obtain the desired identity
for $\gl \in \fapdc$ such that  $\Re \gl + \rho_P$
is $\bar P$-dominant. Now apply analytic continuation.
\qed

\begin{rem}
With the same method of proof
it can be shown
that Lemma \refer{l: comparison Eis} generalizes
to arbitrary parabolic subgroups of $G.$
In the more general lemma, the expression on the left-hand side
is defined as in Harish--Chandra's work, taking account of
Remark \refer{r: gl parameter Eis}. Moreover, the Eisenstein integral
on the right-hand side is defined as in \bib{CDtf}, p.\ 61,
with $\gl$ in place of $-\gl.$
In the proof one has to replace the decomposition
(\refer{e: Iwasawa deco of element})  by
the decomposition induced by $G = K \exp(\fm_P \cap \fp) A_P N_P.$
\end{rem}

\begin{rem}
Lemma \refer{l: comparison Eis} can also be derived from
\bib{BSmult}, Lemma 1, by expressing both Eisenstein integrals as
matrix coefficients of representations of the principal series,
see \bib{BSft}, Eqn.~(25) and \bib{HC1}, Thm.\ 7.1.
\end{rem}

\begin{cor}
\naam{c: sp iso and nE}
Let $P \in \minparabs,$ $\psi \in \cAtwo.$ Then
\begin{equation}
\nE(\spP  \col p^* \psi \col \spgl ) =  p^*( \nE(P \col  \psi \col \gl )) ,
\end{equation}
as an identity of meromorphic functions in the variable $\gl \in \fapdc.$
\end{cor}

\proof
In (\refer{e: comparison Eis}) we substitute
$x_1 = m_1 a$ and $x_2 = m_2 a^{-1}$ for $m_1, m_2 \in M_0 \Ap$ and $a \in \Ap.$
Comparing coefficients
in the asymptotic expansions of type (\refer{e: asymptotics Eis})
for both sides, as $a \to \infty$ in $A_P^+,$ we obtain that
$p^{*-1} \after C_{\spP|\spP}(1 \col \spgl)\after p^* =
C_{P|P}(1 \col \gl)C_{\bar P| P}(1 \col -\bar \gl)^*.$
The result now follows from Lemma \refer{l: comparison Eis} if we apply
the definitions of the normalized Eisenstein integrals,
see (\refer{e: defi normalized Eis}) and  \bib{BSft}, Eqn.\ (49).
\qed

We can now formulate the relation between the Fourier transforms for $G$ and those for the associated
symmetric space $\spspX;$ for the definition of the latter, we refer to \bib{BSft}, Eqn.\ (59).
We define the linear isomorphism
\begin{equation}
\naam{e: p upper star on mero}
p^*: \;\; \Mer(\fapdc) \otimes \cAtwo \;\;\pijl{\simeq}\;\;  \Mer(\spfaqdc) \otimes \spcAtwo
\end{equation}
by
$p^*(\psi)(\spgl) = p^* [\psi(\gl)],$ for $\psi \in \Mer(\fapdc) \otimes \cAtwo $
and for generic $\gl \in \fapdc.$

\begin{lemma}
\naam{l: commutativity diagram Fourier}
Let $P \in \minparabs.$
Then the following diagram commutes:
$$
\begin{array}{ccc}
\Ccminf(G \col \tau) &  \pijl{{\scriptscriptstyle p^*}} & \Ccminf(\spspX \col \sptau)\\
{\scriptscriptstyle\Fou_{\! P}}\downarrow && \downarrow {\scriptscriptstyle\Fou_{\! \spaceind \!P}}\\
\Mer(\fapdc) \otimes \cAtwo & \pijl{{\scriptscriptstyle p^*}}& \Mer(\spfaqdc) \otimes \spcAtwo \; .
\end{array}
$$
\end{lemma}

\proof
Let $f \in \Ccminf(G\col \tau)$ and put $\spf:= p^*(f) \in \Ccminf(\spspX\col \sptau).$
Let $\psi \in \cAtwo.$ Then it follows by application of Corollary \refer{c: sp iso and nE}
and the fact that $dx = p^*(dg)$ that
\begin{eqnarray*}
\inp{\Fou_{\spP}(\spf)(\spgl)}{\sppsi} &=&
\int_{\spspX} \inp{\spf(x)}{\nE(\spP \col \sppsi \col - \spaceind \bar \gl \col x)}\; dx
\\&=&
\int_G \inp{f(g)}{\nE(P\col \psi \col - \bar \gl\col g)}\; dg
\\
&=&
\inp{\Fou_Pf(\gl)}{\psi} = \inp{\spaceind [\Fou_P f (\gl)]}{\sppsi}.
\end{eqnarray*}
In the last equality we have used that $p^*: \psi \mapsto \sppsi$
is a unitary isomorphism from $\cAtwo$ onto $\spcAtwo.$
Using the definition of the map (\refer{e: p upper star on mero})
we conclude that  $\Fou_{[\spP]} \after  p^*(f) =  p^* \after \Fou_P f.$
\qed

The map $p^*: \gl \mapsto \spgl$ is a linear isomorphism
from $\fapd$ onto $\spfaqd,$ mapping the set $\gS$ onto $2 \spgS.$ It follows that the map
$p^*$ in (\refer{e: p upper star on mero})
maps $\Mer(\fapdc, \gS) \otimes \cA_2$  isomorphically onto
$\Mer(\spfaqdc, \spgS) \otimes \spcAtwo.$ Moreover, the transpose of its inverse restricts
to a linear isomorphism
\begin{equation}
\naam{e: sp iso on Lau}
p^*:\;\;
\Mer(\fapdc, \gS)^*_\laur \otimes \cAtwo^* \;\;
\pijl{\simeq}\;\;
\Mer(\spfaqdc, \spgS)^*_\laur \otimes \spcAtwo^*,
\end{equation}
which we shall also denote by $\Lau \mapsto \spLau.$

\begin{lemma}
\naam{l: sp iso of AC}
The isomorphism (\refer{e: sp iso on Lau}) maps $\AC(G, \tau, P)$ onto $\AC(\spspX, \sptau, \spP).$
\end{lemma}

\proof
Let $(x,y) \in G \times G$ and $v \in \Vtau.$ We consider the functions
$
f: \mu \mapsto \dE(\spP\col \mu  \col (x, y))v$ and
$g: \gl \mapsto \dE(P\col \gl \col xy^{-1})v,$
where the first dual Eisenstein integral is defined in a fashion analogous to
(\refer{e: defi normalized dual Eis}), see \bib{BSfi}, Eqn.~(2.3).

It follows from Corollary \refer{c: sp iso and nE} that $f = p^* g.$ Thus, for every
$\Lau \in \Mer(\fapdc)^*_\laur \otimes \cAtwo^*$
we have that
$p^*(\Lau) f = \spLau \spaceind g = \Lau g.$ It follows from this that $\Lau \in \AC(G,\tau, P)$
if and only if
 $p^*(\Lau) \in \AC(\spspX\col \sptau\col \spP).$ The result follows.
\qed

Let $P \in \minparabs$ and $R > 0.$
We define the Paley--Wiener space
 $\PW_R(\spspX\col \sptau\col \spP)$
as in \bib{BSpw}, Def.\ 3.4.
The mentioned definition depends on a choice of positive roots, which we
take to be $\spgS(\spP).$ We enlarge
 this space
to a distributional Paley--Wiener space $\PW^*_R(\spspX\col \sptau\col \spP)$
in complete analogy with
the way in which (b) enlarges (a) in Definition \refer{d: normalized PW}.

\begin{thm}
\naam{t: p star iso between PW}
Let $P \in \minparabs$ and $R > 0.$
The map (\refer{e: p upper star on mero})
restricts to a topological linear isomorphism
\begin{equation}
p^*:\;\;  \PW^{*}_R(G,\tau, P) \;\;\pijl{\simeq}
\;\;\PW^{*}_R(\spspX, \sptau, \spP)
\end{equation}
and to a similar isomorphism between the spaces without the superscript $*.$
\end{thm}

\begin{rem}
\naam{r: real AC}
The definition of $\PW_R(\spspX, \sptau, \spP)$ in \bib{BSpw},
Def.\ 3.4, is not completely analogous to Definition \refer{d: normalized PW}
(a), as the definition in \bib{BSpw} invokes only the relations determined by
the space $\AC_\R(\spspX, \tau, P)$ of Arthur--Campoli functionals with real support;
see also \bib{BSpw}, Def.\ 3.2.
However, it follows by application of \bib{BSpw}, Thm.\ 3.6, that the functions in the
Paley--Wiener
space thus defined satisfy all remaining Arthur--Campoli relations as well.
Consequently, \bib{BSpw}, Def.\ 3.4, determines the same Paley--Wiener space
as the analogue of Definition \refer{d: normalized PW} for the triple $(\spspX, \sptau, \spP).$
A similar remark can be made for the distributional Paley--Wiener space.

It follows from these observations, combined with the results of this paper,
that the Paley--Wiener spaces introduced in Definitions \refer{d: arthurs PW},
\refer{d: unnormalized PW} and \refer{d: normalized PW} remain unaltered if
only the Arthur--Campoli functionals with real support are invoked.
\end{rem}

\proof
We define the hyperplane configuration
$\spHyp = \spHyp_{\spspX, \sptau,\spP}$ and the map
$\spd = d_{\spspX, \sptau, \spP}$ as in \bib{BSpw}, text following
Lemma 2.1.
In view of the relation between the dual
Eisenstein integrals, it follows that
$p^*(\Hyp) = \Hyp$ and that $d_* = d \after p^*.$ This implies that
the map  $p^*$ introduced in  (\refer{e: p upper star on mero})
restricts to a topological linear isomorphism
$$
\cP^*(\fapdc, \Hyp, d) \to \cP^*(\spfaqdc, \spHyp, \spd),
$$
and to a similar isomorphism between the spaces
without the superscript $*.$
In view of Lemma \refer{l: sp iso of AC}
these isomorphisms restrict to isomorphisms of the closed
subspaces with index $\AC.$

Let $\sppi$ be defined as $\pi$ in \bib{BSpw}, for the tuple
$(\spspX, \sptau)$ and the positive
system $\spgS^+ = \spgS(\spP).$ Then it follows from the relation between the dual Eisenstein
integrals that $\sppi = p^*(\pi),$ possibly up to a nonzero constant factor,
which
 we may ignore here.
As $p^*: \fapd \to \spfaqd,$ $\gl \mapsto \spgl,$
is an isometry, it follows that
$$
\pi(\gl) e^{R |\Re \gl|} = \sppi(\spgl) e^{R |\Re \spgl| }.
$$
Moreover, from $p^*(\gS(P)) = 2 \spgS(\spP)$
it follows that $\spfaqd(\spP, 0) = p^*(\fapd(P, 0)).$
Thus, a function $\gf \in \cP_\AC^{(*)}(\fapdc, \Hyp, d)$ satisfies an estimate of
type (\refer{e: PW estimate}) (or of type (\refer{e: PW star estimate})) if
and only if the function $p^*(\gf)$ satisfies the analogous
estimate for the triple $(\spspX, \sptau, \spP).$ The result now follows in view
of Remark \refer{r: real AC}.
\qed

It follows from Lemma  \refer{l: commutativity diagram Fourier} combined with
Lemma \refer{l: sp iso on R supported functions} and
Theorem \refer{t: p star iso between PW} that
$\Fou_P$ is a topological linear isomorphism $C_R^\infty(G\col \tau) \to \PW_R(G,\tau, P)$
if and only if $\Fou_{\spP}$ is a topological linear isomorphism
$C_R^\infty(\spspX\col \sptau) \to \PW_R(\spspX , \sptau,\spP).$
In view of the results
of Section \refer{s: the normalized Fourier transform},
it now follows that
Theorem \refer{t: arthurs PW two}, hence Arthur's Paley--Wiener theorem,
is a consequence of  \bib{BSpw}, Thm.\ 3.6.
Similarly, it follows from the Paley--Wiener theorem proved in \bib{BSdpw} that $\Fou_P$
is a topological linear isomorphism from $C_R^{-\infty}(G\col \tau)$ onto $\PW_R^*(G,\tau, P).$
Thus, the validity of Theorem \refer{t: distributional PW}
follows from the main result of \bib{BSdpw}.

\def\adritem#1{\hbox{\small #1}}
\def\distance{\hbox{\hspace{3.5cm}}}
\def\apetail{@}
\def\adderik{\vbox{
\adritem{E.~P.~van den Ban}
\adritem{Mathematisch Instituut}
\adritem{Universiteit Utrecht}
\adritem{PO Box 80 010}
\adritem{3508 TA Utrecht}
\adritem{Netherlands}
\adritem{E-mail: ban{\apetail}math.uu.nl}
}
}
\def\addhenrik{\vbox{
\adritem{H.~Schlichtkrull}
\adritem{Matematisk Institut}
\adritem{K\o benhavns Universitet}
\adritem{Universitetsparken 5}
\adritem{2100 K\o benhavn \O}
\adritem{Denmark}
\adritem{E-mail: schlicht@math.ku.dk}
}
}
\mbox{\ }
\vfill
\hbox{\vbox{\adderik}\vbox{\distance}\vbox{\addhenrik}}
\end{document}